\def\gtg{\mathfrak{g}}

\def\caU{\mathfrak{U}}
\def\caW{\mathfrak{W}}
\documentclass[12pt]{article}
\usepackage{amsfonts,amssymb,amsmath,fullpage,pst-node,rotating,graphics,epsfig, subfigure}
\usepackage{multirow}
\usepackage{authblk}

\newcommand{\nO}{{\textnormal O}}
\newcommand{\nSO}{{\textnormal {SO}}}
\newcommand{\R}{{\mathbb R}}

\newcommand{\Z}{{\mathbb Z}}

\newcommand{\mH}{{\cal Q}}

\newcommand{\T}{{\mathbb T}}
\newcommand{\I}{{\mathbb I}}
\newcommand{\V}{{\mathbb V}}
\newcommand{\mO}{{\mathbb O}}
\newcommand{\D}{{\mathbb D}}

\newcommand{\cB}{{\cal B}}

\newcommand{\Fix}{{\rm Fix}\,}

\newcommand{\rd}{{\rm d}}

\newcommand{\rl}{{\,|\,}}
\newcommand{\inn}{\textnormal{in}}
\newcommand{\out}{\textnormal{out}}
\newcommand{\abs}[1]{|{#1}|}

\newtheorem{definition}{Definition}
\newtheorem{theorem}{Theorem}
\newtheorem{lemma}{Lemma}

\newtheorem{corollary}{Corollary}
\newtheorem{remark}{Remark}
\newcommand{\proof}{\noindent{\bf Proof: }}

\newcommand{\qed}{\hfill{\bf QED}\vspace{5mm}}
\usepackage{color}
\definecolor{DarkGreen}{rgb}{0.,0.3,0.}

\begin{document}
\title{Simple heteroclinic networks in $\R^4$}

\author[1]{Olga Podvigina}
\author[2]{Alexander Lohse}
\affil[1]{\small Institute of Earthquake Prediction Theory
and Mathematical Geophysics, 84/32 Profsoyuznaya St, 117997 Moscow, Russian Federation}
\affil[2]{\small Universit\"at Hamburg, Fachbereich Mathematik, Bundesstra\ss e 55, 20146 Hamburg, Germany}

\maketitle

\begin{abstract}
We classify simple heteroclinic networks for a $\Gamma$-equivariant system in $\R^4$ with finite $\Gamma \subset \nO(4)$, proceeding as follows: we define a graph associated with a given $\Gamma \subset \nO(n)$ and identify all so-called simple graphs associated with subgroups of $\nO(4)$. Then, knowing the graph associated with a given $\Gamma$, we determine the types of heteroclinic networks that the group admits. Our study is restricted to networks that are maximal in the sense that they have the highest possible number of connections -- any non-maximal network can then be derived by deleting one or more connections. Finally, for networks of type A, i.e., admitted by $\Gamma \subset \nSO(4)$, we give necessary and sufficient conditions for fragmentary and essential asymptotic stability. (For other simple heteroclinic networks the conditions for stability are known.) The results are illustrated by a numerical example of a simple heteroclinic network that involves two subcycles that can be essentially asymptotically stable simultaneously.
\end{abstract}



\section{Introduction}\label{sec1}
Heteroclinic cycles and networks are flow-invariant sets in a dynamical system that are associated with stop-and-go dynamics encountered in a variety of applications. As such they have been studied from various angles over the last decades. The present paper is a systematic investigation of simple robust heteroclinic networks in $\R^4$. Our interest lies in systems of the form
\begin{align}\label{sys1}
\dot{x}=f(x), \text{ with } f(\gamma x)=\gamma f(x) \text{ for all} \ x \in  \R^n,\ \gamma \in \Gamma.
\end{align}
where $f: \R^n \to \R^n$ is a smooth map and $\Gamma \subset \nO(n)$ is a finite group. While some constructions we present are general we mainly focus on the case $n=4$. In a system \eqref{sys1} a heteroclinic cycle is a set of equilibria $\xi_1, \ldots ,\xi_M$ and connections $\kappa_j \subset W^u(\xi_j) \cap W^s(\xi_{j+1})$ in the intersection of the respective unstable and stable manifolds of subsequent equilibria, with the convention $M+1=1$. It is well known that if all connections are of saddle-sink type in fixed-point subspaces $P_j$ of the system, then the cycle persists under $\Gamma$-equivariant perturbations and is therefore said to be robust. A cycle in $\R^4$ is called simple if $\dim P_j=2$ for all $j$, see section \ref{sec2} for details. A (simple) heteroclinic network is a connected union of (simple) heteroclinic cycles.

Simple cycles have previously been studied by many authors, see e.g.\
\cite{cl2000, Kru97} for a general overview. Conditions for their asymptotic
stability were derived in \cite{km95a,km04}. In order to improve our systematic
understanding of heteroclinic dynamics it is desirable to classify
low-dimensional heteroclinic networks, which can possibly be combined to form
more complex structures in higher dimensions. A full classification of simple
cycles in $\R^4$ has been achieved step by step in \cite{km04,pc15, sot03, sot05},
while pseudo-simple cycles (that differ from simple ones in the isotypic
decomposition of $\R^4$ w.r.t.\ certain subgroups) were addressed in
\cite{clp17}. In this paper we take the next step by deriving a complete list of
simple networks in $\R^4$. To do so, we introduce an intuitive way of associating
%a subgroup $\Gamma \subset O(n)$ with a graph: by drawing each isotropy semiaxis
%as a point and each isotropy plane as a line, such that the line contains the
%point if and only if the semiaxis is contained in the respective plane. We say
a subgroup $\Gamma \subset \nO(n)$ with a graph: by drawing
isotropy semiaxes
as points and isotropy planes as lines, such that a line contains a point if and
only if the respective semiaxis is contained in the respective plane. We say
that such a graph is {\it simple} if for any isotropy semiaxis its isotropy
subgroup decomposes $\R^n$ into one-dimensional isotypic components. There is a
close relation between simple graphs associated with $\Gamma \subset \nO(4)$ and
the structure of heteroclinic cycles/networks that a $\Gamma$-equivariant vector
field may exhibit.

In a heteroclinic network no individual cycle can be asymptotically stable, so the situation becomes more subtle than in \cite{km95a,km04}: in the context of networks it is of particular interest to identify subcycles with the strongest possible attraction properties. There are several notions of non-asymptotic stability available in the literature, ranging from the strongest, essential asymptotic stability (e.a.s)\ \cite{Mel91}, to the weakest, fragmentary asymptotic stability (f.a.s.)\ \cite{op12}. We characterise the stability configurations in networks of type A which corresponds to the case $\Gamma \subset \nSO(4)$.

Thus, the two main contributions of this work are (i) a complete description of simple graphs for $\Gamma \subset \nO(4)$, the groups associated with them and the heteroclinic networks they admit, and (ii) a stability study for heteroclinic networks of type A in $\R^4$. The former is achieved by building on results about simple heteroclinic cycles in \cite{pc15}, where the subgroups of $\nO(4)$ have been studied in detail, using the quaternionic representation from \cite{pdv} in order to decide which groups admit heteroclinic cycles. The latter closes a gap in the study of stability in heteroclinic networks in $\R^4$: networks of type Z (as opposed to type A) have been thoroughly investigated before and they give rise to complex stability configurations, \cite{cl14,pa11}. For type A networks we show that the situation is comparatively simple in the sense that fragmentary and essential asymptotic stability of the network are directly linked to the existence of one or more subcycles with the same property.

Moreover, we numerically investigate a heteroclinic network in $\R^4$ consisting of four equilibria and six connections forming two cycles with two equilibria and one cycle with four equilibria. As expected from our theoretic results it is possible to have each cycle e.a.s., but no two cycles sharing a connections can be f.a.s.\ simultaneously.

This paper is organized as follows. In section \ref{sec2} we recall necessary background information and terminology regarding heteroclinic cycles and their stability properties as well as the quaternionic approach for representing subgroups of $\nO(4)$. Section \ref{sec3} explains how we associate a given subgroup $\Gamma\subset \nO(n)$ with a (simple) graph. In Theorems \ref{th32} and \ref{th34} we identify all types of simple graphs that can occur for $\Gamma \subset \nO(4)$ and list the corresponding groups. This enables us to classify maximal simple heteroclinic networks in $\R^4$ in Theorem \ref{th35}. In section \ref{sec4} we give necessary and sufficient conditions for fragmentary and essential asymptotic stability of type A simple heteroclinic networks in $\R^4$. Finally, in section \ref{sec5} we study numerically an example of a heteroclinic network identified in section \ref{sec4} that may exhibit one attracting cycle with four equilibria or two attracting cycles with two equilibria each. Section \ref{sec6} concludes.

\section{Background}\label{sec2}

\subsection{Heteroclinic cycles and notions of stability}
\label{sec:stability}

In this subsection we recall basic terminology in the context of (robust) heteroclinic cycles and their stability properties. Given a $\Gamma$-equivariant dynamical system \eqref{sys1} recall that for $x \in \R^n$ the \emph{isotropy subgroup of $x$} is the subgroup of all elements in $\Gamma$ that fix $x$. On the other hand, for a subgroup $\Sigma \subset \Gamma$ we denote by $\Fix(\Sigma)$ its \emph{fixed point space}, i.e.\ the space of points in $\R^n$ that are fixed by all elements of $\Sigma$.

Let $\xi_1, \ldots ,\xi_M$ be hyperbolic equilibria of a system \eqref{sys1} with stable and unstable manifolds $W^s(\xi_j)$ and $W^u(\xi_j)$, respectively. Also, let $\kappa_j \subset W^u(\xi_j) \cap W^s(\xi_{j+1}) \neq \varnothing$ for $j=1,\ldots,M$ be connections between them, where we set $\xi_{M+1}=\xi_1$. Then the collection of equilibria $\{\xi_1,\ldots ,\xi_M\}$ together with the connecting trajectories $\{ \kappa_1, \ldots ,\kappa_M\}$ is called a \emph{heteroclinic cycle}. A connected union of heteroclinic cycles is a \emph{heteroclinic network}.

A heteroclinic cycle is \emph{structurally stable} or \emph{robust} if for all $j$ there are subgroups $\Sigma_j \subset \Gamma$ such that $\xi_{j+1}$ is a sink in $P_j:=\Fix(\Sigma_j)$ and $\kappa_j$ is contained in $P_j$, see \cite{km95a}. As usual we divide the eigenvalues of the Jacobian $df(\xi_j)$ into \emph{radial} (eigenspace belonging to $P_{j-1} \cap P_j$), \emph{contracting} (belonging to $P_{j-1} \ominus (P_{j-1} \cap P_j)$), \emph{expanding} (belonging to $P_j \ominus (P_{j-1} \cap P_j)$) and \emph{transverse} (all others), where we write $X \ominus Y$ for a complementary subspace of $Y$ in $X$.

We are interested in cycles where (i) $\dim P_j =2$ for all $j$, and (ii) the heteroclinic cycle intersects each connected component of $P_{j-1} \cap P_j=\Fix(\Delta_j)$, where $\Delta_j \subset \Gamma$ is some subgroup, at most once. Note that then $\dim \left(P_{j-1} \cap P_j \right) =1$, and so we refer to $P_j$ as an \emph{isotropy plane} and to $P_{j-1} \cap P_j$ as an \emph{isotropy axis}. The latter is often denoted by $L_j$, we differ slightly from this notation by using $L_j$ for a connected component of $P_{j-1} \cap P_j \setminus \{0\}$, i.e.\ for an \emph{isotropy semiaxis}. Note that different connected components may or may not be related by symmetry.

In $\R^4$, there is then one eigenvalue of each type and we denote the corresponding contracting, expanding and transverse eigenspaces of $df(\xi_j)$ by $V_j$, $W_j$ and $T_j$, respectively. In \cite{pc15} it is shown that under these conditions there are three possibilities for the unique $\Delta_j$-isotypic decomposition of $\R^4$:
\begin{enumerate}
 \item[(1)] $\R^4=\Fix(\Delta_j) \oplus V_j \oplus W_j \oplus T_j$
 \item[(2)] $\R^4=\Fix(\Delta_j) \oplus V_j \oplus \widetilde{W}_j$, where $\widetilde{W}_j=W_j \oplus T_j$ is two-dimensional
 \item[(3)] $\R^4=\Fix(\Delta_j) \oplus W_j \oplus \widetilde{V}_j$, where $\widetilde{V}_j=V_j \oplus T_j$ is two-dimensional
\end{enumerate}
Here $\oplus$ denotes the orthogonal direct sum. This inspires the following definition.

\begin{definition}[\cite{pc15}]
We call a heteroclinic cycle satisfying conditions (i) and (ii) above
\emph{simple} if case 1 holds true for all $j$, and \emph{pseudo-simple} otherwise. A heteroclinic network is called simple if it consists only of simple cycles,
and pseudo-simple if at least one of its cycles is pseudo-simple.
\end{definition}

We aim to identify all subgroups of $\nO(4)$ that admit simple heteroclinic networks in the following sense.

\begin{definition}[\cite{pc15}] \label{def:admits}
We say that a subgroup $\Gamma$ of O($n$) {\em admits} robust heteroclinic
cycles (network) if there exists an open subset of the set of smooth
$\Gamma$-equivariant vector fields in $\R^n$, such that vector fields in this
subset possess a robust heteroclinic cycle (network).
\end{definition}

In order to discuss stability properties a further distinction of simple cycles into different types has proved useful. There are several established ways to do this, we reproduce here only the types that are relevant for our results in section \ref{sec4}.

\begin{definition}[\cite{km04, op12}]
A simple heteroclinic cycle (network) in $\R^4$ is said to be of
 \begin{itemize}
  \item[(i)] type $A$ if $\Sigma_j \cong \Z_2$ for all $j$,
  \item[(ii)] type $Z$ if $\Sigma_j$ decomposes $P_j^\perp$ into one-dimensional isotypic components for all $j$.
 \end{itemize}
\end{definition}
Note that $\Gamma \subset \nSO(4)$ admits heteroclinic cycles only of type A and that
type Z cycles are admitted by $\Gamma\not\subset \nSO(4)$ only. Type A homoclinic
and heteroclinic cycles can exist in $\Gamma$-equivariant systems
with $\Gamma\not\subset \nSO(4)$.

A heteroclinic cycle that belongs to a network cannot be asymptotically stable
because it does not contain the entire unstable manifolds of all its equilibria.
Moreover, it has been proved recently that a compact robust heteroclinic network
comprised of equilibria and a finite number of
connecting trajectories is never asymptotically
stable \cite{PCL2017}.
Various weaker notions of stability have been introduced over the last decades. The strongest one is essential asymptotic stability which goes back to \cite{Mel91}. In the following $\ell(.)$ denotes the Lebesgue measure of a set in $\R^n$ and for $\varepsilon >0$ we denote by $B_{\varepsilon}(X)$ an $\varepsilon$-neighbourhood of $X$. For an $\varepsilon$-neighbourhood of $0$ we simply write $B_\varepsilon$.

\begin{definition} \label{def:eas1}
A set $X$ is {\em essentially asymptotically stable (e.a.s)} if
$$\lim_{\delta\to0}\lim_{\varepsilon\to0}
{\ell(B_{\varepsilon}(X)\setminus{\cal B}_\delta(X))\over
\ell\left(B_{\varepsilon}(X)\right)}=0,$$
where
$$
{\cal B}_\delta(X) = \left\{x\in \R^n~:~ d(\Phi(x,t),X)<\delta\hbox{~for all~}t>0\mbox{~and~} \lim_{t\rightarrow+\infty}d(\Phi(x,t),X)=0 \right\}.
$$
is the $\delta$-basin of attraction of $X$ and $\Phi(x,t)$ denotes the flow of the system.
\end{definition}

Heteroclinic cycles that are not e.a.s.\ may still attract a set of positive Lebesgue measure within their neighbourhood. This is captured in the following term from \cite{op12}.

\begin{definition} \label{def:fragmstable}
A heteroclinic cycle $X$ is {\em fragmentarily asymptotically stable (f.a.s.)} if
$\ell({\cal B}_\delta(X))>0$ for any $\delta>0$.
\end{definition}
We also refer to the (local) stability index along a connection as a way to characterise stability and attraction of heteroclinic cycles and networks. It is a quantity which can be computed with respect to any flow-invariant set and is constant along solution trajectories, first defined in \cite{pa11} and studied further e.g.\ in \cite{Loh15}. In particular, a heteroclinic cycle/network is e.a.s.\ if and only if the local stability indices along all of its connections are positive. It is f.a.s.\ as soon as one of them is greater than $-\infty$.

Finally, we introduce the notions of thin and thick cusps, which we use in our proofs in section \ref{sec4}. For $\alpha>1$ we define the following subset of $\R^2$:
$$V(a_1, a_2, \alpha) := \left\{(x_1,x_2) \in \R^2~:~\abs{a_1x_1+a_2x_2} < \max(\abs{x_1},\abs{x_2})^\alpha \right\}$$
\begin{definition}
We say that $U\subset\R^2$ is a \emph{thin cusp}, if \\
$\bullet$ $\ell(B_{\delta}\cap U)>0$ for all $\delta>0$;\\
$\bullet$ there exist $a_1, a_2$, $\alpha >1$ and $\delta>0$ such that $U \cap B_\delta \subset V(a_1,a_2,\alpha)$.
\end{definition}
\begin{definition}
We say that $U\subset\R^2$ is a \emph{thick cusp}, if its
complement in $\R^2$ is a union of a finite number of thin cusps.
\end{definition}

\begin{remark}\label{rem-cusps}
\begin{enumerate}
 \item Let $U_1, U_2$ be thin cusps. Then, generically, for sufficiently small $\delta >0$, we have $U_1 \cap U_2 \cap B_\delta = \varnothing$.
 \item Let $U_1$ be a thin cusp and $U_2$ be a thick cusp. Then, generically, for sufficiently small $\delta >0$, we have $U_1 \cap B_\delta \subset U_2$.
\end{enumerate}
\end{remark}

\subsection{Quaternions and subgroups of O(4)}
\label{sec:quaternions}

In this section we briefly describe the presentation of finite subgroups
of O(4) with quaternions, for more on this topic see \cite{conw,pdv}.
A real quaternion is a set of four real numbers, ${\bf q}=(q_1,q_2,q_3,q_4)$.
Multiplication of quaternions is defined as
\begin{equation}\label{mqua}
\begin{array}{ccc}
{\bf q}{\bf w}&=&(q_1w_1-q_2w_2-q_3w_3-q_4w_4,q_1w_2+q_2w_1+q_3w_4-q_4w_3,\\
&&q_1w_3-q_2w_4+q_3w_1+q_4w_2,q_1w_4+q_2w_3-q_3w_2+q_4w_1).
\end{array}\end{equation}
The conjugate of $\bf q$ is defined as $\tilde{\bf q}=(q_1,-q_2,-q_3,-q_4)$. For a unit
quaternion we have $\tilde{\bf q}={\bf q}^{-1}$.
We denote by $\mH$ the multiplicative group of unit
quaternions; obviously, its identity element is $(1,0,0,0)$.

Due to the existence of a 2-to-1 homomorphism of $\mH$ onto SO(3),
finite subgroups of $\mH$ are labelled after the respective subgroups of SO(3).
They are:
\begin{equation}\label{finsg}
\renewcommand{\arraystretch}{1.5}
\begin{array}{ccl}
\Z_n&=&\displaystyle{\oplus_{r=0}^{n-1}}(\cos2r\pi/n,0,0,\sin2r\pi/n)\\
\D_n&=&\Z_{2n}\oplus\displaystyle{\oplus_{r=0}^{2n-1}}(0,\cos r\pi/n,\sin r\pi/n,0)\\
\V&=&((\pm1,0,0,0))\\
\T&=&\V\oplus(\pm{1\over2},\pm{1\over2},\pm{1\over2},\pm{1\over2})\\
\mO&=&\T\oplus\sqrt{1\over2}((\pm1,\pm1,0,0))\\
\I&=&\T\oplus{1\over2}((\pm\tau,\pm1,\pm\tau^{-1},0)),
\end{array}\end{equation}
where $\tau=(\sqrt{5}+1)/2$. Double parenthesis denote all even permutations
of quantities within the parenthesis.
Any other finite subgroup of $\mH$ is conjugate to one of these under an inner
automorphism of $\mH$.

For $(q_1,q_2,q_3,q_4)$ regarded as Euclidean coordinates of a point in $\R^4$,
a pair of unit quaternions $({\bf l};{\bf r})$ defines
the transformation ${\bf q}\to{\bf lqr}^{-1}$, which is a rotation in $\R^4$, i.e.
an element of the group SO(4). The mapping
$\Phi:\mH\times\mH\to$\,SO(4) that relates the pair $({\bf l};{\bf r})$
with the rotation ${\bf q}\to{\bf lqr}^{-1}$ is a 2-to-1 homomorphism,
the kernel of which consists of $(1;1)$ and $(-1;-1)$.

Therefore, a finite subgroup of SO(4) is a subgroup of a product of two
finite subgroups of $\mH$.
Denote by $\bf L$ and $\bf R$ the finite subgroups of $\mH$ comprised of
${\bf l}_j$ and ${\bf r}_j$, $1\le j\le J_{\bf l},J_{\bf r}$, respectively.
To any element ${\bf l}\in \bf L$ there are several corresponding elements
${\bf r}_i$, such that $({\bf l};{\bf r}_i)\in\mH$, and similarly for any
${\bf r}\in{\bf R}$. This establishes a correspondence between $\bf L$ and
$\bf R$. Following \cite{pdv}, we denote by ${\bf L}_K$ and ${\bf R}_K$ the
subgroups of $\bf L$ and
$\bf R$ corresponding to the unit elements in $\bf R$ and $\bf L$, respectively,
and write
$({\bf L}\rl{\bf L}_K;{\bf R}\rl{\bf R}_K)$ for the group $\Gamma$.
The isomorphism between ${\bf L}/{\bf L}_K$ and ${\bf R}/{\bf R}_K$ may not be
unique and different isomorphisms may give rise to different subgroups of SO(4).
If this is the case, such subgroups are indicated by additional subscripts
or superscripts.

A reflection in $\R^4$ can be expressed in the quaternionic presentation as
${\bf q}\to{\bf a\tilde qb}$, where ${\bf a}$ and ${\bf b}$ are
unit quaternions. We write this reflection as $({\bf a};{\bf b})^*$.
A group $\Gamma^*\subset$\,O(4), $\Gamma^*\not\subset$\,SO(4), can be decomposed as
$$\Gamma^*=\Gamma\oplus\sigma\Gamma,\hbox{ where $\Gamma\subset$\,SO(4)
and $\sigma=({\bf a};{\bf b})^*\notin$\,SO(4)}.$$

\section{Graphs}\label{sec3}

\subsection{Groups and graphs}\label{sec31}

With a given $\Gamma\subset$\,O($n$) we associate a graph by the following rules:
\begin{itemize}
\item[(A)] A group orbit of isotropy semiaxes is drawn by a point.
\item[(B)] A group orbit of isotropy planes is drawn by a line.
\item[(C)] If for some representatives of the group orbits, an isotropy semiaxis
belongs to an isotropy plane, then the
respective point is drawn on the respective line.
\end{itemize}
\begin{definition}\label{simple}
We say that an isotropy semiaxis of $\Gamma\subset$\,O($n$) is \emph{simple}
if its isotropy group decomposes $\R^n$ into one-dimensional isotypic components.
We call the graph associated with $\Gamma\subset$\,O($n$) \emph{simple} if all
isotropy semiaxes of $\Gamma$ are simple.
\end{definition}
In this subsection we identify all simple graphs associated with finite
subgroups of O(4).

An isotropy semiaxis is a connected component of $\tilde L\setminus\{0\}$, where
$\tilde L$ is an isotropy axis, i.e. $\tilde L=\Fix\Delta$ for some $\Delta\subset\Gamma$.
There are two possibilities for the normalizer $N_{\Gamma}(\Delta)$: either
$N_{\Gamma}(\Delta)=\Delta$ or $N_{\Gamma}(\Delta)=\Delta\times\Z_2$. In the
latter case the two semiaxes comprising $\tilde L\setminus\{0\}$
%are of the same isotropy type,
belong to the same group orbit,
while in the former case they do not.
We use the notion of a semiaxis instead of an axis because in the case
$N_{\Gamma}(\Delta)=\Delta$ an isotropy axis can contain two distinct (i.e.,
not related by a symmetry) steady states, one on each isotropy semiaxis.

The following lemma proves properties of finite groups with simple graphs which
will be used further on to identify graphs associated with subgroups of O(4).
\begin{lemma}\label{lem31}
Consider a finite $\Gamma\subset \nO(n)$, such that
\begin{itemize}
\item[(a)] the graph associated with $\Gamma$ is simple;
\item[(b)] $\Gamma$ has at least one isotropy semiaxis;
\item[(c)] any isotropy semiaxis of $\Gamma$ belongs to an intersection of at least
two isotropy planes.
\end{itemize}
Then
\begin{itemize}
%\item[(i)] Any isotropy plane contains at most two isotropy types of semiaxes.
\item[(i)] Any isotropy plane contains isotropy semiaxes belonging to
at most two group orbits.
\item[(ii)] $\Gamma$ admits simple heteroclinic cycles.
\item[(iii)] If in addition any isotropy semiaxis
% belongs to
is an intersection of at least three isotropy planes and any isotropy plane
%contains exactly two types of isotropy semiaxes,
contains isotropy semiaxes belonging to two group orbits,
then the group admits simple heteroclinic networks.
\item[(iv)] Consider a connected component of the graph associated with
$\Gamma$. If $N_{\Gamma}(\Delta_0)=\Delta_0$ for one isotropy semiaxis of this
component, then $N_{\Gamma}(\Delta_j)=\Delta_j$ for all isotropy semiaxes,
which belong to this component. The component involves at most two
group orbits of semiaxes and each isotropy plane contains
representatives of both orbits.
%isotropy types of semiaxes that belong to all isotropy planes in this component.
\end{itemize}
\end{lemma}

\proof
(i) Let $P_j$ be an isotropy plane of $\Gamma$. Denote by $\Sigma_j$ its
isotropy subgroup and by $N_{\Gamma}(\Sigma_j)$ the
normalizer of $\Sigma_j$ in $\Gamma$. The plane $P_j$ is invariant under
$N_{\Gamma}(\Sigma_j)$ where $N_{\Gamma}(\Sigma_j)/\Sigma_j\cong\D_{K_j}$
for some $K_j\ge0$ (we assume $\D_1\equiv\Z_2$ and $\D_0\equiv I$).
The definition of simple cycles and condition (c) imply that any isotropy semiaxis
that belongs to $P_j$ is an isotropy semiaxis of $N_{\Gamma}(\Sigma_j)/\Sigma_j$.
The group $\D_{K_j}$ splits semiaxes in $P_j$ into two group orbits.
If isotropy semiaxes from different $\D_{K_j}$-group orbits are related by some
$\gamma\in\Gamma$ such that $\gamma\not\in N_{\Gamma}(\Sigma_j)$, then $P_j$
contains semiaxes from just one group orbit.
%The group $\D_{K_j}$ split semiaxes in $P_j$ into two isotropy types.
%If semiaxes of different $\D_{K_j}$-isotropy types are related by some
%$\gamma\in\Gamma$ such that $\gamma\not\in N_{\Gamma}(\Sigma_j)$, then $P_j$
%contains just one type of isotropy semiaxes.

(ii) To begin with, we prove that the group $\Gamma$ has two sequences of
isotropy semiaxes $L_j$ and isotropy planes $P_j$, $1\le j\le m$, such that:
\begin{itemize}
\item[$\bullet$] $L_i\ne \gamma L_j$ and $P_i\ne \gamma P_j$ for any $i\ne j$,
$1\le i,j\le m$ and $\gamma\in\Gamma$.
\item[$\bullet$] $L_j=P_{j-1}\cap P_j$ for $j=1,\ldots,m$, and
$L_1=P_m\cap \gamma P_1$ for some $\gamma\in\Gamma$.
\end{itemize}
(We do not exclude $m=1$, in which case the constructed cycle is homoclinic.)

The existence of such sequences follows from conditions (a) and (b): we take
$L'_1$ to be any isotropy semiaxis of $\Gamma$ and $P'_1$ a plane that $L'_1$
belongs to. We denote by $P'_2$ another isotropy plane that contains $L'_1$
and by $L'_2$ the other isotropy axis in $P'_2$. Repeating this procedure, we
finally obtain a sequence $L'_k,\ldots,L'_{k+m+2}$, such that
$L'_k=\gamma L'_{k+m+2}$ for some $\gamma\in\Gamma$, but $L'_i\ne\sigma L'_j$
for any $\sigma\in\Gamma$ and $i\ne j$, $k\le i,j\le k+m+1$. Denoting
$L_j=L'_{k+j-1}$ and $P_j=P'_{k+j-1}$ we obtain the desired sequences.
The proof of the existence of a $\Gamma$-equivariant system with $\xi_j\in L_j$ and
$\kappa_j\subset P_j$ follows the same ideas as the proof of lemma 1 in
\cite{clp17} and is omitted.

(iii) Decompose the set of isotropy semiaxes into disjoint sets
%${\cal L}_j=\{L_{j1},\ldots,L_{js_j}\}$, where for any $L_{jk}$ and $L_{jl}$
${\cal L}_j=\{\Gamma L_{j1},\ldots,\Gamma L_{js_j}\}$, where for any
$L_{jk}$ and $L_{jl}$
there exist sequences $(\{P_1,\ldots,P_m\};\{L_1,\ldots,L_m\})$ as
constructed in part (ii), where $L_{jk}=L_i$ and $L_{jl}=L_{i'}$ for some
$1\le i,i'\le m$.
%A set can be comprised of just one type of axes.
Here $\Gamma L_{jk}$ is the group orbit of the semiaxis $L_{jk}$.
The set ${\cal L}_j$ can be comprised of just one group orbit.
We assume that the sets ${\cal L}_j$ are maximal, namely that
for any $L_{jk}$ and $L_{j'l}$, where $j\ne j'$, such sequences do not exist.
For a set ${\cal L}_j$ denote by
%${\cal P}_j=\{P_{j1},\ldots,P_{jt_j}\}$ the
${\cal P}_j=\{\Gamma P_{j1},\ldots,\Gamma P_{jt_j}\}$ the
set of isotropy planes, such that $P_{ij}$ contains isotropy
%types
semiaxes from two group orbits, one of which belongs to ${\cal L}_j$ and the other one does not.
The assumption that the sets ${\cal L}_j$ are
maximal implies the existence of at least one ${\cal L}_i$ such that $t_i\le 1$, let
it be ${\cal L}_1$. By the conditions in part (iii) above, ${\cal L}_1$ has
at least two group orbits of isotropy semiaxes and at least one of them
is comprised of axes that do not belong
to $P_{11}$. Denote by $L_1$ an isotropy semiaxis that does not belong to $P_{11}$.
Since $L_1$ is an intersection of at least three isotropy planes, different
from $P_{11}$, there exist sequences
$(\{P_1,\ldots,P_m\};\{L_1,\ldots,L_m\})$
and $(\{P'_1,\ldots,P'_{m'}\};\{L_1,\ldots,L_{m'}\})$ obtained by the same
procedure as in part (ii) with $L_j\in{\cal L}_1$ and $L'_j\in{\cal L}_1$,
such that $L_i\ne L_j$ for any $2\le i\le s$ and $2\le j\le s'$, and
$L_{s+k}=L'_{s'+k}$ for any $1\le k\le m-s=m'-s'$.
Similar to \cite{clp17} and part (ii), we construct a dynamical system with
$\xi_j\in L_j$ and $\xi'_i\in L'_i$ (where $\xi_j=\xi'_i$ for $j> s$)
and $\kappa_j\in P_j$ and $\kappa'_i\in P'_i$ (where $\kappa_j=\kappa'_i$
for $j> s$).

(iv) If $L_1$ and $L_m$ belong to the same connected component of the graph,
then there exist sequences $(\{P_1,\ldots,P_{m-1}\};\{L_1,\ldots,L_m\})$ such
that $L_j\subset P_j$ and $L_{j+1}\subset P_j$. If
$N_{\Gamma}(\Delta_1)=\Delta_1$ then $K_j$ is odd and $-I\not\in \D_{K_1}$,
which implies that
$N_{\Gamma}(\Delta_2)=\Delta_2$. Hence, $-I\not\in \D_{K_2}$, which implies that
$N_{\Gamma}(\Delta_3)=\Delta_3$. Repeating this procedure $m-1$ times, we
obtain that $N_{\Gamma}(\Delta_m)=\Delta_m$.

Since all $K_j$ are odd, any $P_j$ in the above sequence contains just one
isotropy type of axis, $\tilde L_j$, and two connected components of
$\tilde L_j\setminus\{0\}$
are not related by symmetries of $\Gamma$.
%of different isotropy types.
Hence, all $\tilde L_j$, $1\le j\le m$, are of the
same isotropy type and the connected components of $\tilde L_j\setminus\{0\}$
for all $1\le j\le m$
% are of two isotropy types.
belong to two distinct group orbits.
\qed

\begin{theorem}\label{th32}
Let $\Gamma$ be a finite subgroup of SO(4), such that $\Gamma$ has at least one
isotropy axis and all its isotropy axes are simple. Then the graph associated
with $\Gamma$ is non-empty and simple and it is one of those shown in figure 1.
The subgroups of SO(4) with associated simple graphs and the types of these
graphs are listed in (\ref{listth1}) below.

\hskip -2cm
%\begin{table}[h]
\begin{equation}\label{listth1}
\renewcommand{\arraystretch}{1.5}
\begin{array}{|l|c|}
\hline
\hbox{Group}&\hbox{Graph}\\
\hline
(\D_{2K_1}\rl\D_{2K_1};\D_{2K_2}\rl\D_{2K_2}), K_1\wedge K_2=1,
K_1+K_2\hbox{ even}& \hbox{V}\\
(\D_{2K_1}\rl\D_{2K_1};\D_{2K_2}\rl\D_{2K_2}), K_1\wedge K_2=1,
K_1+K_2\hbox{ odd}& \hbox{VI}\\
(\D_{2K_1r}\rl\Z_{4K_1};\D_{2K_2r}\rl\Z_{4K_2})_s,\
K_1,K_2,r,s\hbox{ satisfy (\ref{condno1})}& \hbox{IV}\\
(\D_{2K_1r}\rl\Z_{2K_1};\D_{2K_2r}\rl\Z_{2K_2})_s,\ K_1,K_2\hbox{ odd, }
K_1,K_2,r,s\hbox{ satisfy (\ref{condno2})}&\hbox{III}\\
(\D_{2K_1}\rl\D_{K_1};\D_{2K_2}\rl\D_{K_2}), K_1\wedge K_2=1&\hbox{IV}\\
(\D_{2K_1}\rl\D_{K_1};\D_{2K_2}\rl\Z_{4K_2}),\
K_1\hbox{ even, }K_1/2\wedge K_2=1 &\hbox{II}\\
(\D_{2K_1}\rl\D_{K_1};\D_{2K_2}\rl\Z_{4K_2}),\ K_1\hbox{ odd},\
K_1\wedge K_2=1 &\hbox{III}\\
(\D_{2K}\rl\D_{2K};\T\rl\T),\ K\hbox{ even}&\hbox{I}\\
(\D_{2K}\rl\D_{2K};\T\rl\T),\ K\hbox{ odd}&\hbox{II}\\
(\D_{2K}\rl\D_{2K};\mO\rl\mO),\ K\hbox{ odd},\ K\ne3k &\hbox{III}\\
(\D_{2K}\rl\Z_{4K};\mO\rl\T),\ K\ne3k,\ K\hbox{ even}&\hbox{I}\\
(\D_{2K}\rl\Z_{4K};\mO\rl\T),\ K\ne3k,\ K\hbox{ odd}&\hbox{II}\\
(\D_{2K}\rl\D_K;\mO\rl\T),\ K\ne2(2k+1),\ K\ne3k,\ K\hbox{ even}&\hbox{I}\\
(\D_{2K}\rl\D_K;\mO\rl\T),\ K\ne2(2k+1),\ K\ne3k,\ K\hbox{ odd}&\hbox{II}\\
(\D_{2K}\rl\D_{2K};\I\rl\I),\ K\ne5k,\ K\hbox{ even}&\hbox{I}\\
(\D_{2K}\rl\D_{2K};\I\rl\I),\ K\ne5k,\ K\hbox{ odd}&\hbox{II}\\
(\D_{2K_1r}\rl\Z_{K_1};\D_{2K_2r}\rl\Z_{K_2})_s,\
K_1,K_2\hbox{ odd,}K_1,K_2,r,s\hbox{ satisfy (\ref{condno3})}&\hbox{II}\\
\hline
\end{array}
\end{equation}
The conditions mentioned are
\begin{equation}\label{condno1}
K_1\wedge K_2=1,\quad r\wedge K_2-sK_1=1,
\end{equation}
\begin{equation}\label{condno2}
K_1\wedge K_2=1,\quad r\wedge (K_2\pm sK_1)/2=1,
\end{equation}
\begin{equation}\label{condno3}
K_1\wedge K_2=1,\quad r\wedge (K_2\pm sK_1)/2=1,\quad r\wedge (K_2\pm sK_1)/4=1,
\end{equation}
where plus or minus are taken so that the ratios are integer and $K_1 \wedge K_2$ denotes the greatest common divisor of $K_1$ and $K_2$.
%\end{table}
\end{theorem}

\begin{figure}[h]
%[p]

\vspace*{-2mm}
\hspace*{-10mm}\includegraphics[width=6cm]{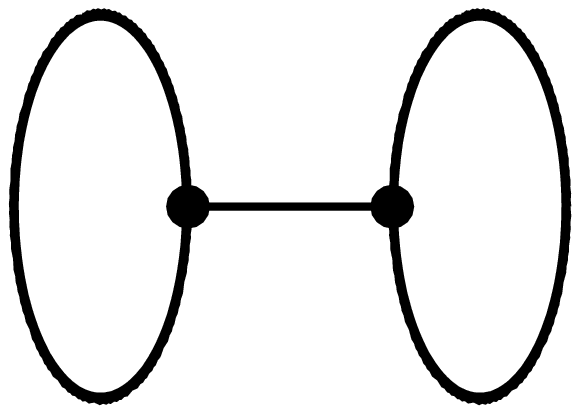}\hspace*{-3mm}
\includegraphics[width=6cm]{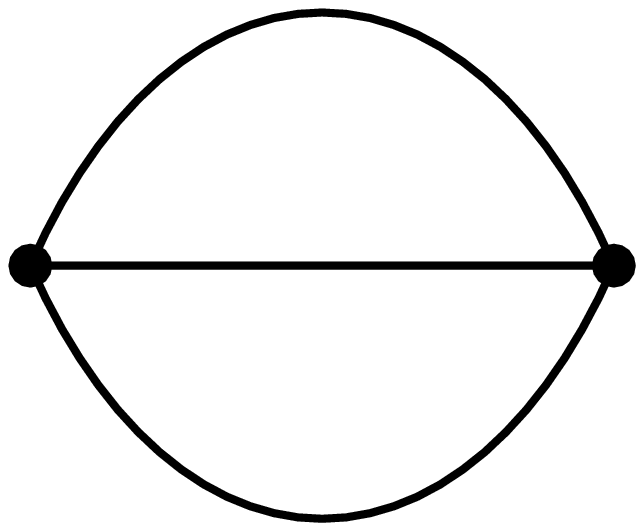}\hspace*{-3mm}
\includegraphics[width=6cm]{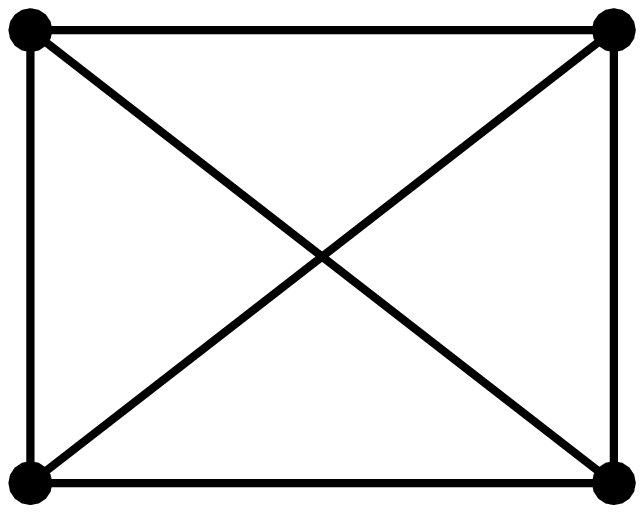}

\vspace*{-12mm}
\hspace*{40mm}{\Large I}\hspace*{53mm}{\Large II}\hspace*{50mm}{\Large III}

\vspace*{5mm}
\hspace*{-13mm}\hspace*{5mm}\includegraphics[width=6cm]{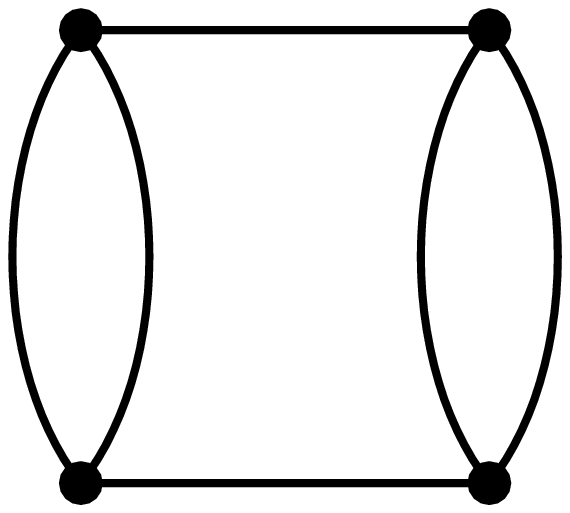}\hspace*{-3mm}
\includegraphics[width=6cm]{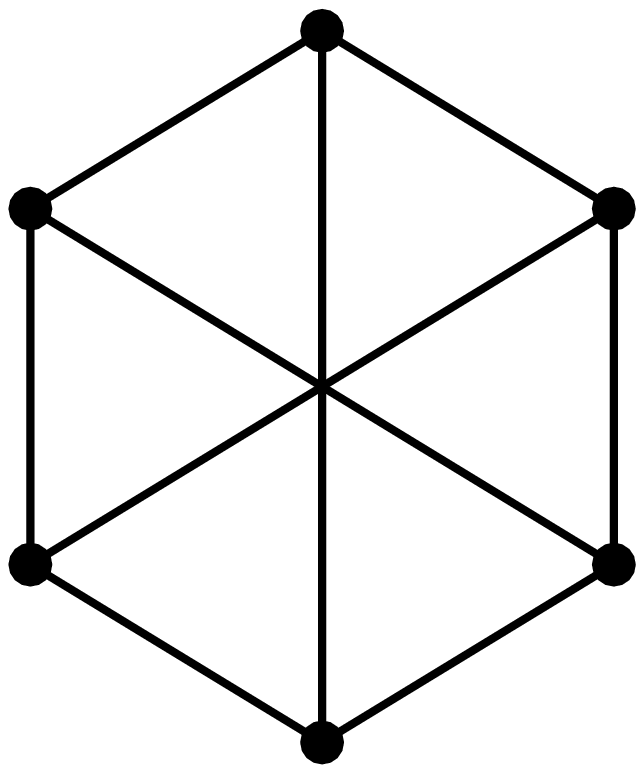}\hspace*{-3mm}\includegraphics[width=6cm]{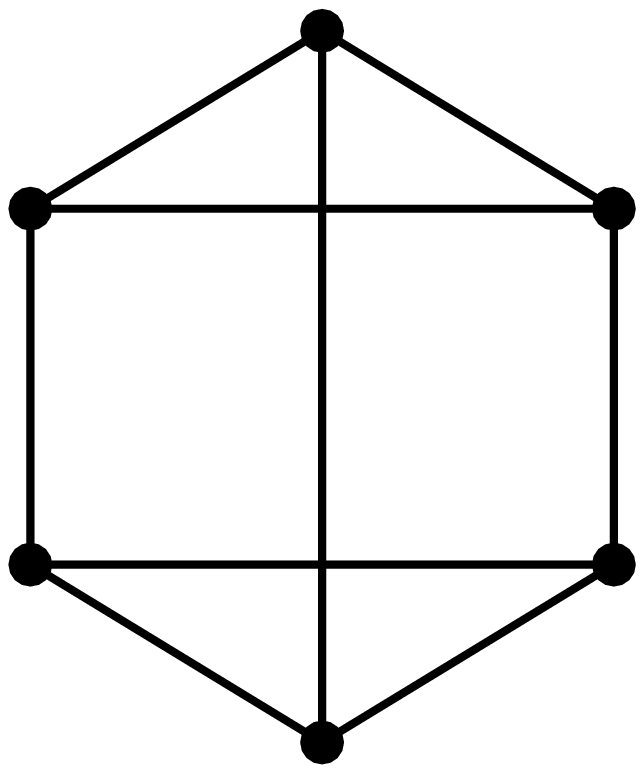}

\vspace*{-10mm}
\hspace*{40mm}{\Large IV}\hspace*{53mm}{\Large V}\hspace*{50mm}{\Large VI}

\vspace*{4mm}
\noindent
\caption{Simple graphs associated with subgroups of O(4).}
\label{fig1}\end{figure}

\proof
The graph associated with $\Gamma$ is simple and non-empty, because all
isotropy axes are simple and it has at least one isotropy axis. Since any plane has
% at most two types of isotropy axes,
isotropy semiaxes belonging to at most two group orbits,
we draw the points representing the semiaxes
as endpoints of the line representing the plane. So, instead of points we can
use the term {\it vertices}. A line that corresponds
to an isotropy plane either connects two vertices, or it begins and ends at the
same vertex. Since $\Gamma$ is a subgroup of SO(4), a simple isotropy axis
belongs to an intersection of three isotropy planes. Hence, the number of
isotropy types of semiaxes is 2/3 of the number of isotropy types of planes and
the number of isotropy semiaxes is even.

Our proof uses results of \cite{pc15}, where simple heteroclinic cycles in
$\R^4$ were studied. The table in appendix C {\it ibid} lists conjugacy classes of
isotropy subgroups of finite groups $\Gamma\subset$\,SO(4) satisfying
$\dim\Fix\,(\Sigma)=2$ and $\dim\Fix\,(\Delta)=1$. Such groups $\Sigma$ always
satisfy $\Sigma\cong\Z_2$. In the table only selected $\Delta$ such that
$\Delta\cong(\Z_2)^2$ are listed.

Since $P=\Fix\,(\Sigma)$ implies that
$\gamma P=\Fix\,(\gamma\Sigma\gamma^{-1})$, the number of conjugacy classes
of $\Sigma$ equals the number of group orbits of isotropy planes.
%From the table we see that the number of isotropy types of intersecting planes
From the table we see that the number of group orbits
of intersecting isotropy planes
for various $\Gamma\subset$\,SO(4) can be 3, 6 or 9. Hence the number of
%isotropy types of
group orbits of
semiaxes is 2, 4 or 6. (If an isotropy plane does not
intersect with another isotropy plane, then it does not contain a heteroclinic
connection and, therefore, can be ignored.)
If the number of
% isotropy types of
group orbits of
semiaxes is 2, an isotropy plane can be
homoclinic (i.e.,
% type of isotropy semiaxis).
 all its isotropy semiaxes belong to the same group orbit).
In such a case there exists another
(not related by a symmetry of $\Gamma$)
%(isotropy type of)
homoclinic plane
and the graph is of type I. Another possibility for a graph
with two vertices is type II where two vertices are connected by three lines.

For a connected graph with four vertices and without homoclinic isotropy planes
(connectedness and absence of homoclinic planes follow from appendix C
in \cite{pc15}) there are two possibilities: a graph without two-vertices
connected subgraphs (type III) and the one that involves such subgraphs
(type IV). There can possibly exist several types of graphs with six
vertices, however, as we conclude from table C, the only groups that have six
isotropy types of axes are $(\D_{2K_1}\rl\D_{2K_1};\D_{2K_2}\rl\D_{2K_2})$. Explicitly constructing
graphs associated with these groups using data from tables in appendices B and
C in \cite{pc15}, we obtain that
the graph is of type V when $K_1+K_2$ is even and of type VI when it is odd.

Subgroups of SO(4) admitting simple heteroclinic cycles are listed
in theorem 2 in \cite{pc15}. For consistency of presentation we reproduce
this list below:
\hskip -2cm
%\begin{table}[h]
%\begin{equation}\label{listold}
$$
\renewcommand{\arraystretch}{1.5}
\begin{array}{|l|l|l|}
\hline
(\D_{2K_1}\rl\D_{2K_1};\D_{2K_2}\rl\D_{2K_2}) &
(\D_{2K}\rl\Z_{4K};\mO\rl\T),\ K\ne3k \\
\hline
(\D_{2K_1r}\rl\Z_{4K_1};\D_{2K_2r}\rl\Z_{4K_2})_s,\ K_1,K_2,r,s\hbox{ satisfy (\ref{condno1})}&
(\D_{2K}\rl\D_K;\mO\rl\T),\ K\ne2(2k+1)\\
\hline
(\D_{2K_1r}\rl\Z_{2K_1};\D_{2K_2r}\rl\Z_{2K_2})_s,\  K_1\wedge K_2=1,&
(\D_{2K}\rl\D_{2K};\I\rl\I),\ K\ne5k\\
\cline{2-2}
K_1,K_2,r,s\hbox{ satisfy (\ref{condno2})}&
(\T\rl\Z_2;\T\rl\Z_2)\\
\hline
(\D_{2K_1}\rl\D_{K_1};\D_{2K_2}\rl\D_{K_2})&
(\T\rl\T;\mO\rl\mO) \\
\hline
(\D_{2K_1}\rl\D_{K_1};\D_{2K_2}\rl\Z_{4K_2}),\ K_1\hbox{ even, }
K_1/2\wedge K_2=1 &
(\mO\rl\mO;\I\rl\I) \\
\hline
(\D_{2K_1}\rl\D_{K_1};\D_{2K_2}\rl\Z_{4K_2}),\ K_1\hbox{ odd}&
(\D_{2K_1r}\rl\Z_{K_1};\D_{2K_2r}\rl\Z_{K_2})_s,\
K_1,K_2\hbox{ odd,}\\
\cline{1-1}
(\D_{2K}\rl\D_{2K};\T\rl\T)&
K_1,K_2,r,s\hbox{ satisfy (\ref{condno3})}\\
\hline
(\D_{2K}\rl\D_{2K};\mO\rl\mO),\ K\hbox{ odd} &\\
\hline
\end{array}
$$
%\end{equation}
%\end{table}

In this list some subgroups (e.g., $(\T\rl\T;\mO\rl\mO)$ or
$(\D_{2K_1}\rl\D_{2K_1};\D_{2K_2}\rl\D_{2K_2})$ where $K_1$ and $K_2$ are not
co-prime) have non-simple isotropy axes. For
$\Gamma\subset$\,SO(4) any non-simple axis is an intersection of two planes,
that are not orthogonal. The existence of isotropy planes, intersecting non-orthogonally,
can be identified using appendices B and C and lemma 3 in \cite{pc15}.
Excluding from the above table subgroups that have isotropy planes with
non-orthogonal intersections we obtain the list (\ref{listth1}) of subgroups that have
associated simple graphs.

We use the data from appendices B-D in \cite{pc15} to identify the graph for each
group. The groups with three
% types of isotropy planes we split
group orbits of isotropy planes we split
into types I and II depending on whether the group admits homoclinic cycles or
not. The groups with six
%types  of isotropy planes split into types
group orbits of isotropy planes split into types
of III and IV, depending on whether they have two-vertices subgraphs or not.
For the groups $(\D_{2K_1}\rl\D_{2K_1};\D_{2K_2}\rl\D_{2K_2})$, the only ones with nine isotropy types of
planes, we construct graphs from the data in the tables, as stated in the
beginning of the proof of this lemma.
\qed

\begin{theorem}\label{th34}
The graph associated with a group $\Gamma^*\subset$\,O(4),
$$\Gamma^*=\Gamma\oplus\sigma\Gamma,
\hbox{ where $\Gamma\subset$\,SO(4) and $\sigma\notin$\,SO(4)},$$
is non-empty and simple if and only if ${\Gamma}$ and $\sigma$ are listed in (\ref{listth2}).
\begin{table}[h]
\begin{equation}\label{listth2}
\renewcommand{\arraystretch}{1.5}
\begin{array}{|l|l|l|}
\hline
{\Gamma} & \sigma \\
\hline
(\D_2\rl\Z_2;\D_2\rl\Z_2) & ((0,1,0,0),(0,1,0,0))^* \\
\hline
(\D_2\rl\Z_1;\D_2\rl\Z_1) & ((1,0,0,0),(1,0,0,0))^* \\
\hline
\end{array}
\end{equation}
\end{table}
\end{theorem}

\proof
If the graph associated with a group $\Gamma^*$ is non-empty and simple,
then its subgroup $\Gamma$ admits simple heteroclinic cycles. According to
theorem 3 in \cite{pc15}, such $\Gamma$ is one of the following:
\begin{equation*}
\renewcommand{\arraystretch}{1.5}
\begin{array}{c}
(\D_4\rl\Z_2;\D_4\rl\Z_2),\ (\D_4\rl\Z_1;\D_4\rl\Z_1)_3,\
(\D_2\rl\Z_2;\D_2\rl\Z_2),\\
(\T\rl\Z_2;\T\rl\Z_2),\ (\D_2\rl\Z_1;\D_2\rl\Z_1),\
(\D_{2K}\rl\D_K;\D_{2K}\rl\D_K).
\end{array}
\end{equation*}
In the above list, the only groups where all isotropy axes are simple
are the two given in (\ref{listth2}). The respective elements $\sigma$
are also given in \cite{pc15}. The graphs, associated with $\Gamma^*$, are
the same as the ones associated with $\Gamma$, namely type III for
$(\D_2\rl\Z_2;\D_2\rl\Z_2)$ and type II for $(\D_2\rl\Z_1;\D_2\rl\Z_1)$.
\qed

Note that not every simple heteroclinic cycle is admitted by a group $\Gamma$ that is associated with a simple graph. This is because $\Gamma$ may admit simple heteroclinic cycles, but have additional non-simple isotropy axes. However, it can be deduced from tables B-D in \cite{pc15} that all such groups admit at most one simple heteroclinic cycle, and thus no simple heteroclinic network. Therefore, our method allows us to find all simple heteroclinic networks in $\R^4$.

\subsection{Graphs and networks}\label{sec32}

In this subsection we address the following question: For a given
$\Gamma\subset$\,O(4), what kinds of heteroclinic networks can a
$\Gamma$-equivariant dynamical systems possess? For the group $\Gamma^*=\Gamma\oplus\sigma\Gamma\cong \Z_2^4$, with $\Gamma=(\D_2\rl\Z_2;\D_2\rl\Z_2)$ and the corresponding $\sigma$ from theorem \ref{th34}, this question was addressed in \cite{cl16b}, where it was shown that the group supports three different types of networks associated with type III: a network of two cycles with three equilibria each (and one common connection), a network of one cycle with three and one with four equilibria (and two common connections), and a network of two cycles with three and one with four equilibria. In the latter network every available isotropy plane contains a heteroclinic connection, while in the former two there is always a plane that is not required for the construction.

Heteroclinic cycles and networks can be represented by diagrams similar to the graphs
in figure \ref{fig1}. In these diagrams a steady state is drawn as a point, a heteroclinic connection a a line and arrows indicate the directions of the heteroclinic connections. Although our graphs are formally different from
the diagrams (a diagram is drawn for a dynamical system that possesses a
heteroclinic cycle or network, while a graph is drawn for a subgroup of O($n$)),
the graphs can be related to the heteroclinic diagrams using arguments
similar to the ones employed to prove lemma 1 in \cite{clp17}. Namely, we
can explicitly construct a $\Gamma$-equivariant dynamical system such that each isotropy semiaxis contains a steady state and each isotropy plane contains a heteroclinic connection, which can go in either direction.
Hence, we can relate heteroclinic cycles to the graphs shown in figure \ref{fig1}. Note that the geometric structure of the graph, associated
with a group, determines the geometric structure of the networks admitted by this group.
In this subsection we consider maximal networks in the sense that {\bf C1:}
\begin{itemize}
\item[$\bullet$] any isotropy semiaxis contains a steady state,
\item[$\bullet$] any isotropy plane contains a heteroclinic connection,
\item[$\bullet$] all such steady states and connections belong to the network.
\end{itemize}
The group $\Gamma$ can admit non-maximal networks, see e.g., \cite{cl14}
or \cite{cl16b}, where the graph of type III was considered, however all other
networks can be constructed by removing one or more connections from a maximal network.

\begin{figure}[h]
%[p]

\vspace*{-2mm}
\hspace*{13mm}\includegraphics[width=5cm]{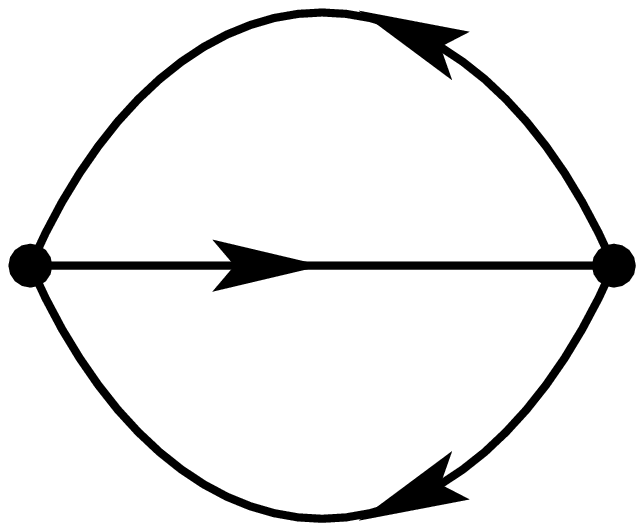}\hspace*{13mm}
\includegraphics[width=5cm]{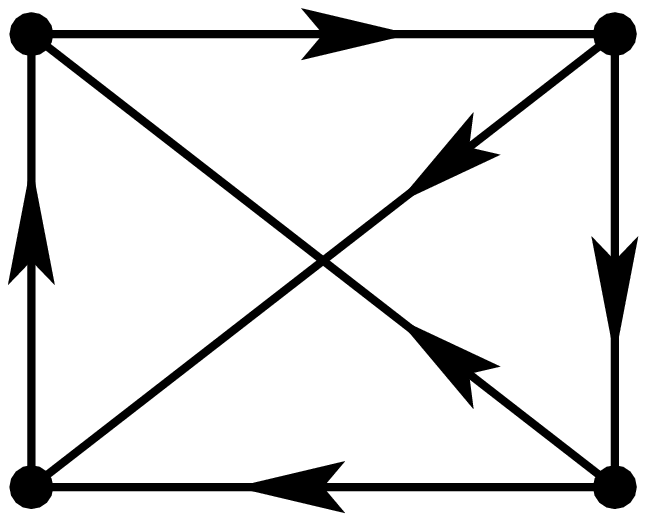}

\vspace*{-11mm}
\hspace*{54mm}{\Large II}\hspace*{55mm}{\Large III}

\vspace*{2mm}
\noindent
\caption{Heteroclinic networks of types II and III.}
\label{fig2}\end{figure}

\begin{figure}[h]
%[p]

\vspace*{-7mm}
\hspace*{-3mm}\includegraphics[width=5cm]{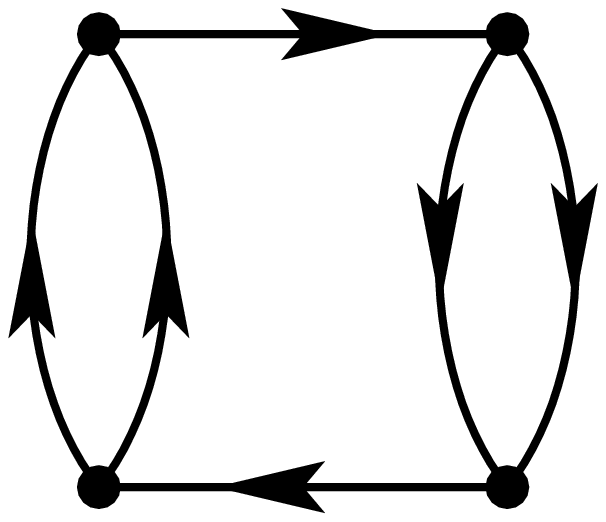}\hspace*{1mm}
\includegraphics[width=5cm]{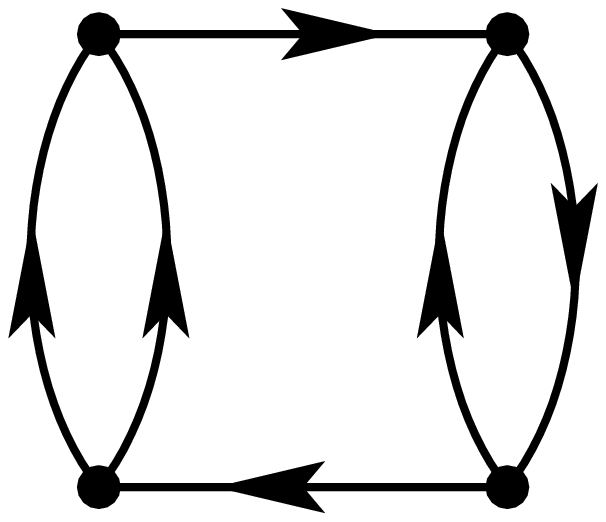}
\hspace*{1mm}\includegraphics[width=5cm]{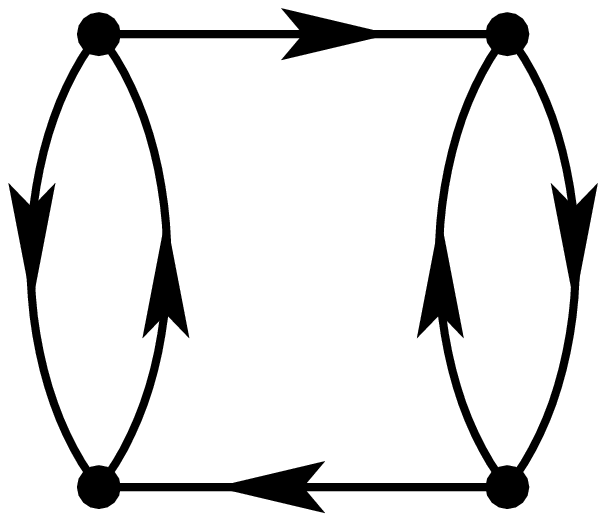}

\vspace*{-11mm}
\hspace*{28mm}{\Large IVa}\hspace*{43mm}{\Large IVb}\hspace*{43mm}{\Large IVc}

\vspace*{2mm}
\noindent
\caption{Heteroclinic networks of type IV.}
\label{fig3}\end{figure}

\begin{figure}[h]
%[p]

\vspace*{-3mm}
\hspace*{10mm}\includegraphics[width=4.5cm]{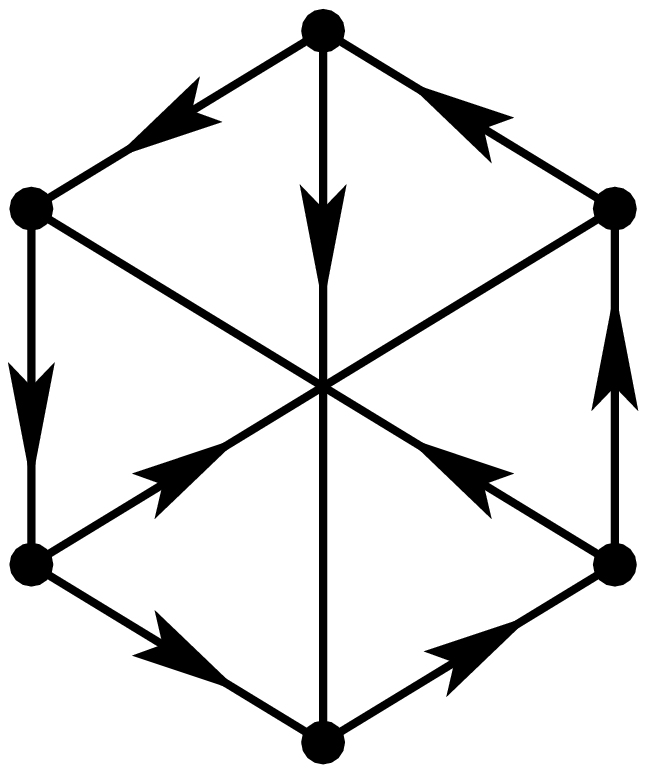}\hspace*{5mm}
\includegraphics[width=4.5cm]{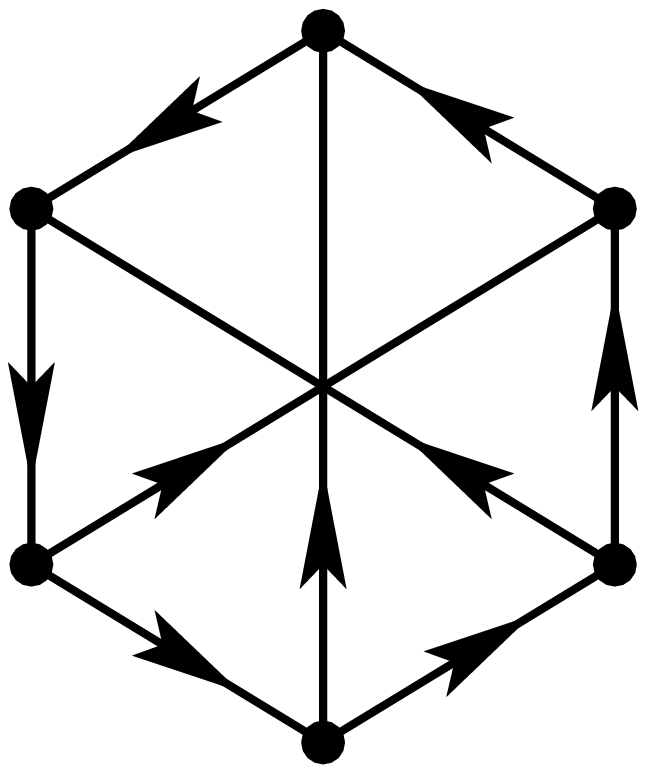}
\hspace*{5mm}\includegraphics[width=4.5cm]{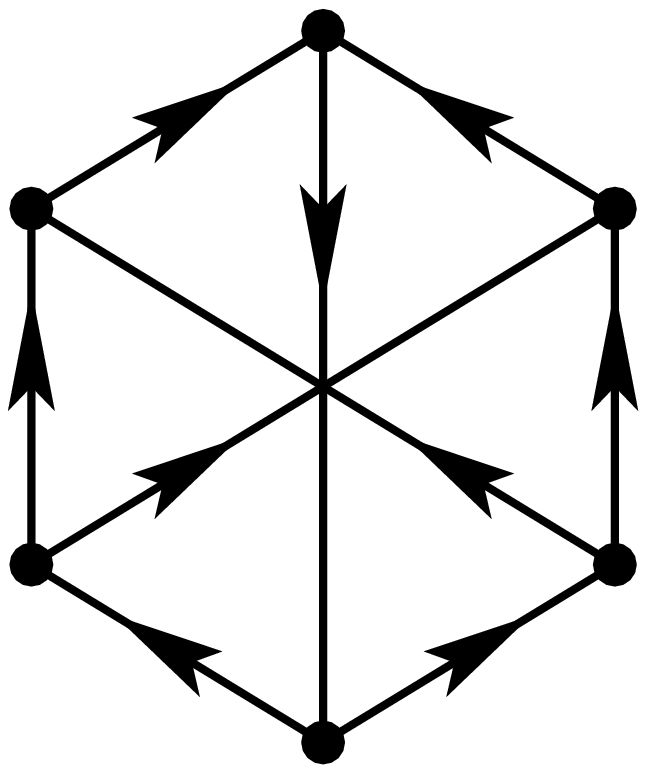}

\vspace*{-10mm}
\hspace*{47mm}{\Large Va}\hspace*{45mm}{\Large Vb}\hspace*{45mm}{\Large Vc}

\vspace*{2mm}
\noindent
\caption{Heteroclinic networks of type V.}
\label{fig4}\end{figure}

\begin{figure}[h]
%[p]

\vspace*{3mm}
\hspace*{-15mm}\includegraphics[width=4.5cm]{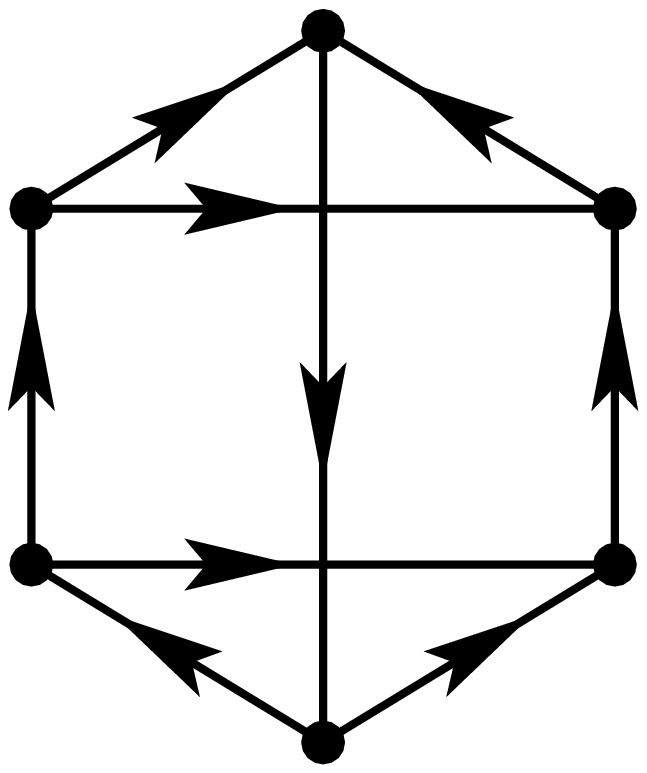}\hspace*{2mm}
\includegraphics[width=4.5cm]{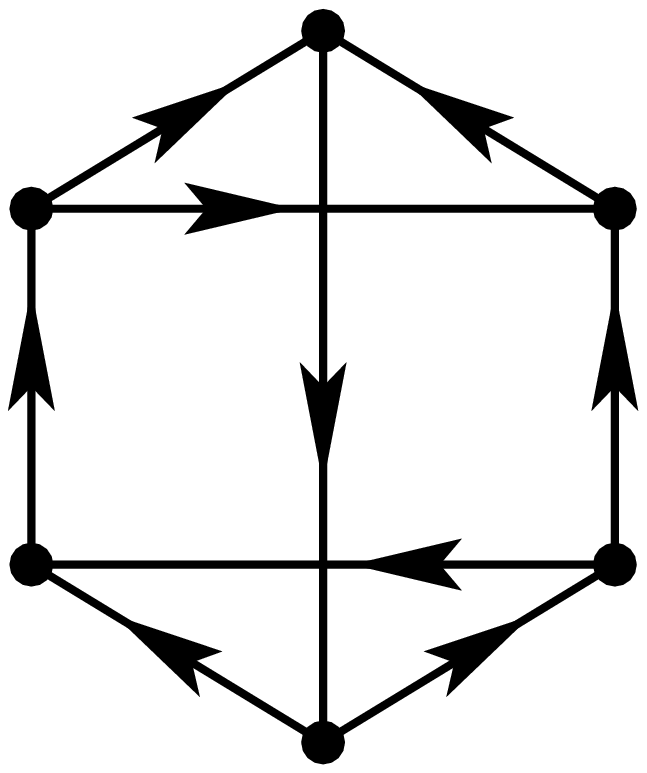}
\hspace*{2mm}\includegraphics[width=4.5cm]{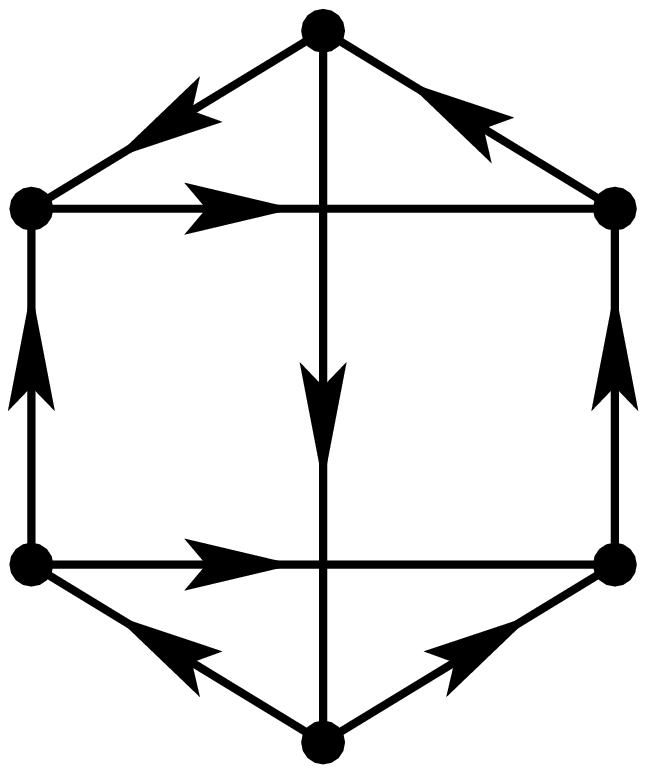}\hspace*{2mm}\includegraphics[width=4.5cm]{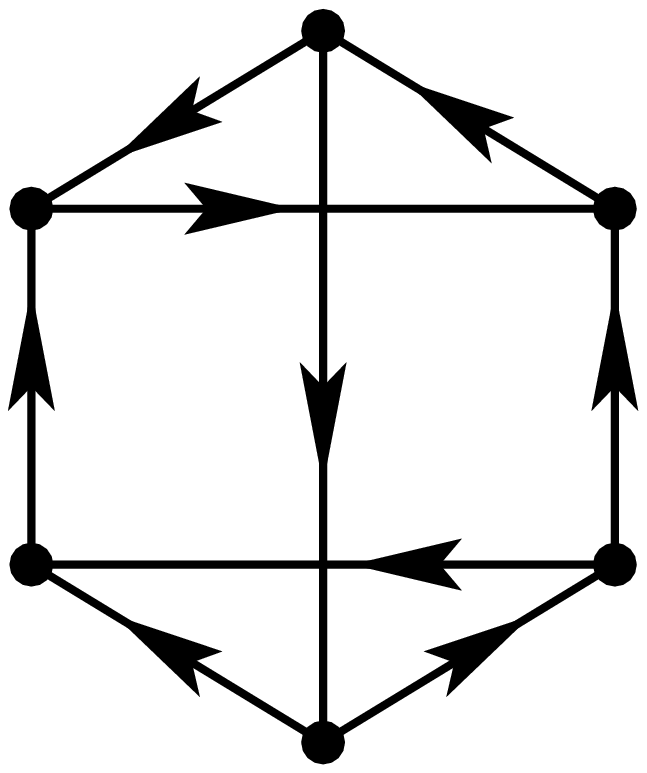}

\vspace*{-9mm}
\hspace*{24mm}{\Large VIa}\hspace*{39mm}{\Large VIb}
\hspace*{39mm}{\Large VIc}\hspace*{39mm}{\Large VId}

\vspace*{5mm}
\hspace*{-15mm}
\includegraphics[width=4.5cm]{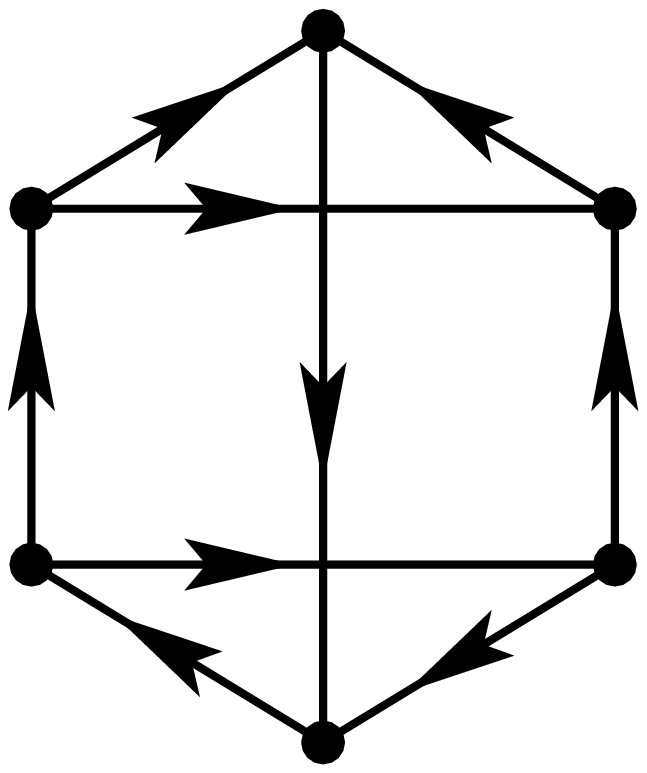}
\hspace*{2mm}\includegraphics[width=4.5cm]{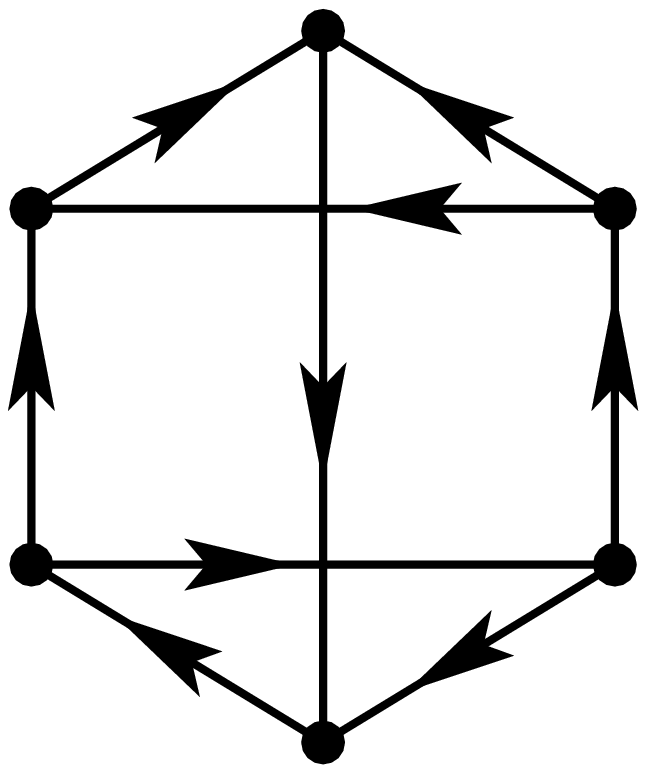}
\hspace*{2mm}\includegraphics[width=4.5cm]{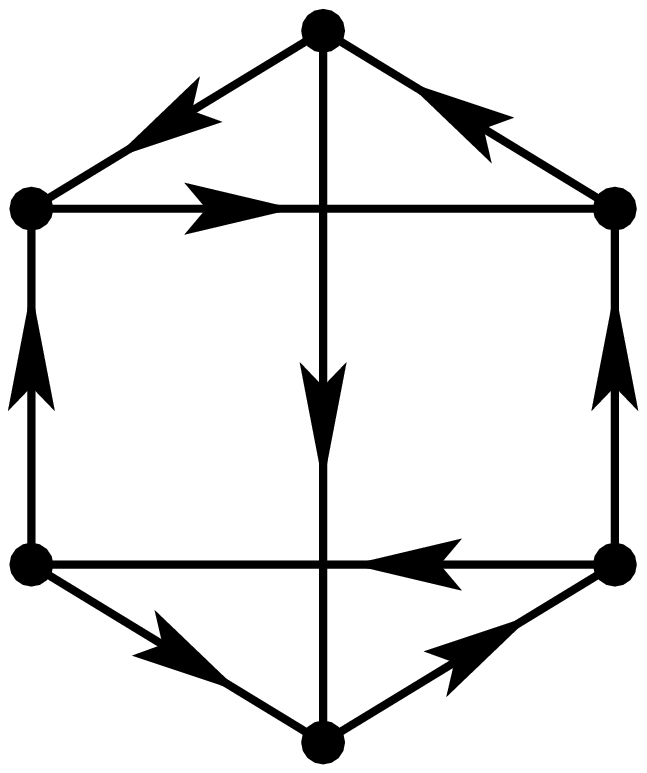}\hspace*{2mm}
\includegraphics[width=4.5cm]{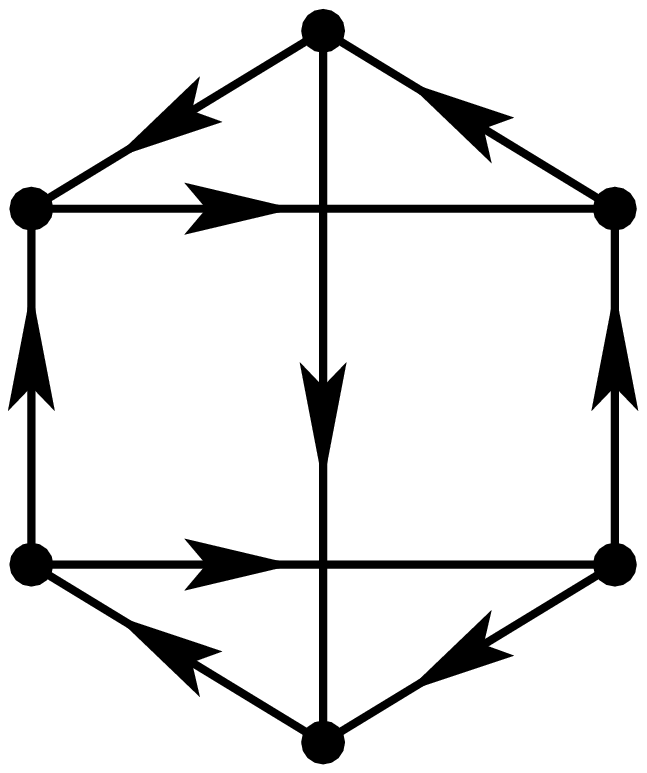}

\vspace*{-9mm}
\hspace*{24mm}{\Large VIe}\hspace*{39mm}{\Large VIf}
\hspace*{39mm}{\Large VIg}\hspace*{39mm}{\Large VIj}

\vspace*{2mm}
\noindent
\caption{Heteroclinic networks of type VI.}
\label{fig5}\end{figure}

\begin{theorem}\label{th35}
Suppose that a $\Gamma$-equivariant system, where $\Gamma$ is a finite
subgroup of SO(4), possesses a heteroclinic network that is simple and
maximal. Then
\begin{itemize}
\item[(i)] If $\Gamma$ has exactly two
group orbits of isotropy semiaxes, then the network is
%types of isotropy axes, then the network is
of type II in figure \ref{fig2}.
\item[(ii)] If the graph associated with $\Gamma$ is of type III, then the
 network is of type III in figure \ref{fig2}.
\item[(iii)] If the graph associated with $\Gamma$ is of type IV, then the
 network is one of types IVa-IVc in figure \ref{fig3}.
\item[(iv)] If the graph associated with $\Gamma$ is of type V, then the
 network is one of types Va-Vc in figure \ref{fig4}.
\item[(v)] If the graph associated with $\Gamma$ is of type VI, then the
 network is one of types VIa-VIg in figure \ref{fig5}.
\end{itemize}
\end{theorem}

\proof
(i) If $\Gamma$ has two
group orbits of isotropy semiaxes, then the associated graph
%types of isotropy axes, then the associated graph
is of type I or II. The groups with associated graphs of type I
admit two distinct (i.e., not related by a symmetry of $\Gamma$) homoclinic
cycles, however they do not admit heteroclinic networks.
The graph of type II involves two
group orbits of semiaxes,
%isotropy types of semiaxes,
and three
group orbits of planes, hence, the heteroclinic networks admitted by
%types of planes, hence, the heteroclinic networks admitted by
the respective groups involve two
(group orbits of) equilibria, $\xi_1$ and
%isotropy types of steady states, $\xi_1$ and
$\xi_2$ (one for each of isotropy semiaxes), and three
(group orbits of) heteroclinic connections
(one for each of isotropy planes), say, two from $\xi_1$ to $\xi_2$ and one from
$\xi_2$ to $\xi_1$.

(ii) If the graph associated with $\Gamma$ is of type III, then the network
is a union of graphs involving three or four equilibria.
%isotropy types of equilibria.
Evidently, it has at least one cycle with three equilibria,
% isotropy types
let it be $\xi_1\to\xi_2\to\xi_3\to\xi'_1$. There are four possibilities for
the remaining connections:
% (isotropy types of) connections:
$$(a)\ \xi_3\to\xi_4,\ \xi'_4\to\xi_1,\ \xi''_4\to\xi_2;\quad
(b)\ \xi_3\to\xi_4,\ \xi'_4\to\xi_1,\ \xi_2\to\xi''_4;$$
$$(c)\ \xi_3\to\xi_4,\ \xi_1\to\xi'_4,\ \xi''_4\to\xi_2;\quad
(d)\ \xi_3\to\xi_4,\ \xi'_4\to\xi_1,\ \xi_2\to\xi''_4.$$
The network (b) is mapped to the network (a) by the permutation of
equilibria $(\xi_1,\xi_2,\xi_3,\xi_4)\mapsto(\xi_4,\xi_1,\xi_2,\xi_3)$;
the network (c) is mapped to the network (a) by the permutation of
equilibria $(\xi_1,\xi_2,\xi_3,\xi_4)\mapsto(\xi_4,\xi_2,\xi_3,\xi_1)$;
the network (d) is mapped to the network (a) by the permutation of
equilibria $(\xi_1,\xi_2,\xi_3,\xi_4)\mapsto(\xi_3,\xi_1,\xi_2,\xi_4)$.
Hence, the only possibility is the network (a) that is shown in figure \ref{fig2}.

(iii) If the graph associated with $\Gamma$ is of type IV, then the network
can involve zero, one or two cycles that are comprised of two
% isotropy types of
equilibria. Hence, there are three different heteroclinic networks, shown in figure \ref{fig3}.

(iv) If the graph is of type V, the network either does or does not involve a cycle with
six equilibria. In case there is such a cycle, say $\xi_1\to\xi_2\to\ldots\to\xi_6\to\xi'_1$,
there are two possibilities (up to cyclic permutation of equilibria)
for the remaining connections: the equilibria unstable in the transverse directions
are $(\xi_1,\xi_2,\xi_3)$ or $(\xi_1,\xi_3,\xi_5)$.
If the network does not involve a cycle with six equilibria, there is only
one such network up to a permutation of $\xi_j$. The proof is
similar to that of part (ii) and we omit it. The latter network can be
thought of as a union of four cycles, each connecting four equilibria. These networks
are shown in figure \ref{fig4}.

(v) The graph of type VI can be considered as a union of two subgraphs (upper and lower) with three vertices each, and three vertical lines connecting these
subgraphs. Let the equilibria of the upper subgraph be $(\xi_1,\xi_2,\xi_3)$
and those of the lower $(\xi_4,\xi_5,\xi_6)$ such that vertical lines connect
$\xi_j$ and $\xi_{j+3}$. According to the definition of a network,
two of the vertical connections must go in the same direction and the third
one in the opposite. We assume there are connections $\xi_1\to\xi_4$,
$\xi_2\to\xi_5$ and $\xi_6\to\xi_3$. Then there are four possibilities for the
connections between $\xi_1$, $\xi_2$ and $\xi_3$:
$$(a)\ \xi_1\to\xi_2,\ \xi_2\to\xi_3,\ \xi_3\to\xi'_1;\quad
(b)\ \xi_1\to\xi_2,\ \xi_3\to\xi_2,\ \xi_3\to\xi'_1;$$
$$(c)\ \xi_2\to\xi_1,\ \xi_3\to\xi_2,\ \xi_3\to\xi'_;\quad
(d)\ \xi_2\to\xi_1,\ \xi_3\to\xi_2,\ \xi'_1\to\xi_3.$$
Similarly, there are four possibilities for the connections on the lower
subgraph. The only permutation of $\xi_j$ that preserves the directions of
connections on the vertical lines is
$(\xi_1,\xi_2,\xi_3,\xi_4,\xi_5,\xi_6)\mapsto(\xi_2,\xi_1,\xi_3,\xi_5,\xi_4,\xi_6)$,
hence the number of different networks is eight.
They are shown in figure \ref{fig5}.
\qed

\section{Stability of type A heteroclinic networks}
\label{sec4}

In this section we consider an equivariant dynamical system \eqref{sys1} in $\R^4$ that possesses a simple heteroclinic network of type A. We give necessary and sufficient conditions for its fragmentary and essential asymptotic stability below. In accordance with the terminology in \cite{ac98} and \cite{km95a}, we use the following definitions.

\begin{definition} \label{def:dom}
Let $Y$ be a simple heteroclinic network in $\R^n$ and $\xi_j\in Y$ an
equilibrium in $Y$. The eigenvalue of $\rd f(\xi_j)$ with maximal real part is called the
\emph{principal eigenvalue at $\xi_j$}. The corresponding eigenspace is the \emph{principal eigenspace}, and the branch of $W^u(\xi_j)$ that is tangent to it is the \emph{principal unstable manifold}, denoted by $W^{pu}(\xi_j)$. A cross section to $W^{pu}(\xi_j)$ near $\xi_j$ is called \emph{principal outgoing cross section}. If $W^{pu}(\xi_j)$ connects to another equilibrium in $Y$, we speak of a \emph{principal connection}. Finally, a cycle comprised only of principal connections is called a \emph{principal cycle}.
\end{definition}

\begin{definition} \label{def:m-p}
Let $Y$ be a simple heteroclinic network in $\R^n$ and $\xi_j\in Y$ an
equilibrium in $Y$. The eigenvalue of $\rd f(\xi_j)$ with minimal real part, excluding the
radial ones, is called the \emph{minus-principal eigenvalue at $\xi_j$}. As in the previous definition we extend this terminology to \emph{minus-principal eigenspaces, cross sections} and \emph{cycles}.
\end{definition}

\begin{remark}\label{rem41}
With this terminology, recall the following two facts from \cite{pa11}.
\begin{itemize}
 \item[(i)] If a type A heteroclinic cycle in $\R^4$ is f.a.s., then it is e.a.s.
 \item[(ii)] If a type A heteroclinic cycle in $\R^4$ is f.a.s., then it is principal.
\end{itemize}
\end{remark}

In the following, let $Y =(\cup_{1\le j\le J} \xi_j)\cup(\cup_{(ij)\in Q}\kappa_{ij}) \subset \R^4$
be a simple heteroclinic network, consisting of equilibria $\xi_j$ and
connections $\kappa_ {ij}$. We use the standard notation for cross sections
$H_{ij}^\inn$ and $H_{jk}^\out$ to incoming and outgoing connections
$\kappa_{ij}$ and $\kappa_{jk}$ near $\xi_j$.
These can be considered as two-dimensional since the radial directions are
irrelevant in the study of stability properties.
We denote the principal outgoing cross section near $\xi_j$ by $H^{p}_j$.

For equilibria $\xi_i, \xi_j, \xi_k \in Y$ with connections $\kappa_{ij}, \kappa_{jk} \subset Y$ we employ the usual notation for local maps $\phi_{ijk}:H_{ij}^\inn \to H_{jk}^\out$, global maps $\psi_{jk}: H_{jk}^\out \to H_{jk}^\inn$ and their compositions $g_{ijk}:=\psi_{jk }\circ \phi_{ijk}: H_{ij}^\inn \to H_{jk}^\inn$. Note that this is a slight abuse of notation because the domain of $\phi_{ijk}$ is only a subset of $H_{ij}^\inn$ if there is more than one positive eigenvalue at $\xi_j$. If there is no ambiguity we drop subscripts and write $\phi_j, \psi_j$ and $g_j$ or even simply $\phi, \psi$ and $g$.

\begin{lemma}\label{lem2-new}
Consider a local map $\phi: H_{ij}^{\inn} \to H_j^{p}$. Then, generically, there
exists a thick cusp $V \subset H_{ij}^{\inn}$ such that
\begin{itemize}
 \item[(i)] $\phi(V)$ is a thin cusp;
 \item[(ii)] for any $\delta>0$ we can find $\varepsilon>0$ such that
$\phi(B_{\varepsilon}\cap V)\subset B_{\delta}$.
\end{itemize}
\end{lemma}
\proof
The local map is given by
$(y_1,y_2)=\phi(x_1,x_2)=(Ax_1^\alpha,x_2|x_1|^\beta)$, where $A \neq 0$, $\alpha >0$ and $\beta >-1$, see e.g.\ \cite{pa11}. Suppose first that $\alpha <1 +\beta$. Set
$$V:=\{(x_1,x_2) \in \R^2~:~\abs{x_2} < \abs{x_1}^h \ \text{and} \ \abs{x_2}<\abs{x_1}^{-\beta+s} \},$$
where $h<1$, $s>0$, $\alpha<h+\beta$ and $-\beta+s<1$. Then $V$ is a thick cusp and for $(y_1,y_2) \in \phi(V)$ we obtain
$$\abs{y_2}<\abs{x_1}^h\abs{x_1}^\beta = A'\abs{y_1}^{\frac{h+\beta}{\alpha}}$$
for suitable $A'>0$, which shows that $\phi(V)$ is a thin cusp.

Now suppose that $\alpha >1+\beta$. In this case we set
$$V:=\{(x_1,x_2) \in \R^2~:~\abs{x_2} > \abs{x_1}^h \ \text{and} \ \abs{x_2}<\abs{x_1}^{-\beta+s} \},$$
where $h>1$, $s>0$, $\alpha>h+\beta$ and $-\beta+s<1$. Again $V$ is a thick cusp and for $(y_1,y_2) \in \phi(V)$ we obtain
$$\abs{y_2}>\abs{x_1}^h\abs{x_1}^\beta = A'\abs{y_1}^{\frac{h+\beta}{\alpha}}$$
for suitable $A'>0$, which as before means that $\phi(V)$ is a thin cusp.
Hence, (i) is proven. Note that $s>0$ is chosen such that $V$ is contained in
the domain of $\phi$ if $\beta<0$. To prove (ii), we note that for
$(x_1,x_2)\in V$ we have $|\phi(x_1,x_2)|<k(\max(|x_1|,|x_2|))^{\min(\alpha,s)}$
with a constant $k>0$.
\qed

\begin{corollary}\label{principal}
Suppose that a set $U\subset H^{\inn}_{ij}$ is a thin cusp and $\kappa_{jk}$
is a principal connection.
\begin{itemize}
 \item[(i)] Then, generically, for any small $\delta>0$ there exists $\varepsilon>0$ such that for all points $x\in B_{\varepsilon}\cap U \setminus \{0\}$ the trajectories $\Phi(x,t)$ leave the $\delta$-neighbourhood of $\xi_j$ through $H^{p}_j=H_{jk}^\out$.
 \item[(ii)] The set $g_j(B_\varepsilon \cap U) \subset H_{jk}^\inn$ is generically a thin cusp.  We also express this as ``thin cusps generically follow principal connections''.
\end{itemize}
\end{corollary}
\proof (i) Since $\kappa_{jk}$ is principal, the domain of $\phi_j$ is a thick cusp in $H_{ij}^\inn$. Then since $U$ is a thin cusp, by remark \ref{rem-cusps} the set $U \cap B_\varepsilon$ is generically contained in it for sufficiently small $\varepsilon >0$.

(ii) By lemma \ref{lem2-new} and (i) above the set $\phi_j(B_\varepsilon \cap U) \subset H_{jk}^\out$ is a thin cusp. The global map $\psi_j$ is a generic linear map, so it maps a thin cusp in $H_{jk}^\out$ into a thin cusp in $H_{jk}^\inn$.
\qed

\begin{corollary}\label{cor4}
Let $X$ be an f.a.s.\ type A heteroclinic cycle in $\R^4$ with equilibria $\xi_k$, $1 \leq k \leq K$. Then,
generically, for all $k$ and $\delta >0$ the set $\cB_{\delta}(X)\cap H_{k-1,k}^\inn$ is a thick cusp.
\end{corollary}
\proof
By lemma \ref{lem2-new} there is a thick cusp $V \subset H_{k-1,k}^\inn$ such
that its image $\phi_k(V) \subset H_{k,k+1}^\out$ is a thin cusp. Then by
corollary \ref{principal} (ii) $\phi_k(V)$ generically follows principal
connections. Since $X$ is f.a.s., the trajectories passing through $V$ stay near $X$ for all
positive times. Therefore, $V \cap B_\varepsilon \subset \cB_{\delta}(X)\cap H_{k-1,k}^\inn$ for sufficiently small $\varepsilon >0$.
\qed

We are now in a position to prove necessary and sufficient conditions for essential asymptotic stability of $Y$. Recall that there are many ways in which type Z networks can be e.a.s., see \cite{cl14,pa11}. For type A networks the situation is simpler as the following theorem shows.

\begin{theorem}\label{th42}
Generically, a type A heteroclinic network $Y$ in $\R^4$ is e.a.s.\ if and only if the following holds.
\begin{itemize}
 \item[(i)] All principal subcycles of $Y$ are f.a.s.
 \item[(ii)] For any equilibrium $\xi_j\in Y$, the network contains a principal connection $\kappa_{jk} \subset Y$.
\end{itemize}
\end{theorem}

\proof
Assume first that (i) and (ii) are satisfied.
Then by remark \ref{rem41} every principal subcycle of $Y$ is e.a.s.
For any connection $\kappa_{12}$ there is a unique sequence of equilibria $\xi_3, \ldots \xi_k \in Y$ with
principal connections $\kappa_{j,j+1}$ for $j=2, \ldots ,k-1$, such that
$\xi_l$ belongs to a principal subcycle $X\subset Y$ if and only if $l= k$.
By lemma \ref{lem2-new} there is a thick cusp $V \subset H_{12}^\inn$ such that $\phi_2(V) \subset H_{23}^\out \equiv H_2^p$ is a thin cusp. By corollary \ref{principal}, thin cusps generically follow principal connections, so $g_k \circ \ldots \circ g_2(V) \subset H_k^\out$ is also a thin cusp. By remark \ref{rem-cusps}, for sufficiently small $\varepsilon >0$, the set $g_k \circ \ldots \circ g_2(V) \cap B_\varepsilon$ is generically contained in the thick cusp that is $\cB_\delta(X) \cap H_k^\out$ (corollary \ref{cor4}). Applying these
arguments to all connections in $Y$, we obtain that $Y$ is e.a.s.

We now prove the other implication. Let $Y$ be e.a.s\ , $X=\xi_1\to...\to\xi_K\subset Y$ be a principal cycle and
$V\subset H_{1}^\inn$ be the set considered in lemma \ref{lem2-new}.
Since $V$ is a thick cusp and $Y$ is e.a.s., the intersection
$\cB_\delta(Y)\cap V$ has positive measure. By corollary \ref{principal} almost all trajectories starting there stay near
$X$ for all $t>0$. Hence, $X$ is f.a.s., which proves (i). If (ii) is not
satisfied, then there exists $\xi_j\in Y$, such that its principal unstable
manifold does not belong to the network. Then by lemma \ref{lem2-new} there is
a thick cusp in $H_{ij}^\inn$ that is mapped to the principal outgoing cross
section and thus does not stay near $Y$, which contradicts $Y$ being e.a.s.
\qed

We now turn to fragmentary asymptotic stability. Given $\delta>0$
we denote by $\caW$ the set of positive measure subsets
$U\subset \cB_{\delta}(Y)\cap H_{ij}^\inn$. In what follows we drop the superscript ``in'' to simplify notation.
By $\gtg$ we define the map acting on subsets $\caU\subset\caW$,
$\caU=\{U_1,...,U_K\}$, as follows:
\begin{itemize}
\item[$\bullet$] If $\caU=\{U\}$, where $U\subset H_{ij}$, and $\xi_j$ has
$M$ outgoing connections, $\kappa_{jj_1},...,\kappa_{jj_M}$, then
$\gtg\caU=\{U_1,...,U_K\}$ and
\begin{equation}
\renewcommand{\arraystretch}{1.5}
\begin{array}{ll}
U_k=&\{~x\in H_{jj_m}\ :\ \Phi(y,0)\in U,\quad\Phi(y,t_0)=x
\hbox{ for some }t_0>0\\
&~~\hbox{ and }\Phi(y,t)\notin H_{i'j'}
\hbox{ for any }0<t<t_0\hbox{ and }(i'j')\in Q~\}$$
\end{array}
\end{equation}
Since the sets $U_k$ have positive measures, we have $0<K\le M$.
\item[$\bullet$] $\gtg(\caU_1\cup\caU_2)=(\gtg\caU_1)\cup(\gtg\caU_2)$
\end{itemize}

\begin{lemma}\label{lem4n1}
Let $U\subset H_{ij}\cap \cB_{\delta}(Y)$ be a thin cusp. Then generically
for $\caU=\{B_\varepsilon\cap U\}$, any $l>0$
and sufficiently small $\varepsilon >0$ we have $\gtg^l\caU=\{U_l\}$,
where $U_l$ is a thin cusp.
\end{lemma}

\proof
We start with the case $l=1$.
The equilibrium $\xi_j$ can have one or two outgoing connections.
First, we assume that there is one outgoing connection, $\kappa_{j,j+1}$.
The map $g:H_{ij}\to H_{j,j+1}$ that maps $x\in H_{ij}$ to
$\Phi(x,t)\in H_{j,j+1}$ is known -- see, e.g., \cite{pa11} and also our comments before lemma \ref{lem2-new}. It is
$$g(x_1,x_2)=(A_1 x_1^{\alpha}+A_2 x_2|x_1|^{\beta},A_3 x_1^{\alpha}+A_4 x_2|x_1|^{\beta}),$$
where $\alpha>0$, $\alpha$  and $\beta$ depend on eigenvalues of $\rd f(\xi_j)$ and
generically $\alpha\ne\beta+1$, $A_j\ne 0$ and $A_1/A_2\ne A_3/A_4$.

Assume first that $\alpha<\beta+1$.
By definition of a thin cusp, for any $c>0$ and small $\varepsilon=\varepsilon(c)>0$
the set $U$ satisfies $B_\varepsilon\cap U\subset W_c$, where
$$W_c=\{~(x_1,x_2)\in\R^2~:~|x_1+\frac{a_2}{a_1}x_2|<c|x_2|~\}.$$
We have
$$g(qx_2,x_2)=(A_1 q^{\alpha}x_2^{\alpha}+A_2 q^{\beta}x_2^{1+\beta},
A_3 q^{\alpha}x_2^{\alpha}+A_4 q^{\beta}x_2^{1+\beta}).$$
Therefore, for any $0<s<(1+\beta)/\alpha-1$ and sufficiently small $x_2$,
any $(y_1,y_2) \in gW_c \subset H_{j,j+1}$ satisfies
$$|y_1-{A_1\over A_3}y_2|<|y_1|^{(1+\beta)/\alpha-s},$$
which implies that $gW_c$ is a thin cusp. Hence,
$U_1=g(B_\varepsilon\cap U)$ is a thin cusp. The arguments in the case $\alpha>\beta+1$ are similar.

In the case of two outgoing connections from $\xi_j$, for small $\delta$ the
trajectories from $U$ (except for a set of zero measure) follow
the principal connection of $\xi_j$. Hence,
$\gtg U=\{ g(B_\varepsilon\cap U)\}$, where $g(B_\varepsilon\cap U)$
is a thin cusp as it is shown above and $g$ is the map along the principal connection.

Applying these arguments $l$ times we prove the lemma.
\qed

\begin{lemma}\label{lem4n2}
Suppose that $\xi_j\in X$ has two incoming connections, $\kappa_{j_1j}$
and $\kappa_{j_2j}$. Let $U^1\subset H_{j_1j}$ and $U^2\subset H_{j_2j}$ be such that $\ell(B_\varepsilon\cap U^k)>0$ for any $\varepsilon>0$ and $k=1,2$.
Then generically at least one of the sets $W^1$ and $W^2$, where
$\{W^k\}=\gtg\{U^k\}$, is a thin cusp.
\end{lemma}

\proof
The maps $g_k:H_{j_kj}\to H_{ji}$ are
$$g_k(x_1,x_2)=(A_{k1}x_1^{\alpha_k}+A_{k2}x_2|x_1|^{\beta_k},A_{k3}x_1^{\alpha_k}+A_{k4}x_2|x_1|^{\beta_k}),$$
where $\beta_k>0$, $\beta_1=\alpha_2$, $\beta_2=\alpha_1$ and generically
$\alpha_1\ne\beta_1$ and $A_j\ne0$. Hence, we have either $\alpha_1<\beta_1$ or
$\alpha_2<\beta_2$. Suppose that $\alpha_1<\beta_1$. Then for $s>0$ such
that $\alpha_1<\beta_1-s$ and sufficiently small $\varepsilon$ the points
$(y_1,y_2)\in g_1(U^1\cap B_\varepsilon)$ satisfy
$$|y_1-{A_{11}\over A_{13}}y_2|<|y_2|^{(\beta_1-s)/\alpha_1},$$
which implies that $W^1=g_1U^1$ is a thin cusp.
\qed

\begin{lemma}\label{lem4n3}
Let $X_1,...,X_L$ be minus-principal cycles that are subsets of $Y$.
(If there are no such cycles, $L=0$ is assumed.)
Consider a heteroclinic connection $\kappa_{ij}\subset Y$ such
that $\kappa_{ij}\not\subset X_l$ for any $1\le l\le L$.
Then there exists $\delta>0$ and $M>0$ such that for $\caU=\{ U\}$,
where $U=\cB_{\delta}(Y)\cap H_{ij}$ we have
$$\gtg^M\caU=\{U_{M1},...,U_{MK}\},$$
where all $U_{Mk}$, $1\le k\le K$, are thin cusps.
\end{lemma}

\proof
Suppose the statement of the theorem does not hold true.
We take $M$ to be the number of all connection comprising $Y$
plus one and consider
$$\gtg^M\caU=\{U_{M1},...,U_{MK}\}.$$
For $\delta\to 0$ the set remains non-empty, hence there exists some
$U_{Mk}$ that is not a thin cusp for any small $\delta$.
The set $U_{Mk}$ is the image of a subset of $U$ under $g=g_M...g_1$, where
$g_m:H_{j_{m-1}j_m}\to H_{j_mj_{m+1}}$, with $j_0=i$ and $j_1=j$, is the map
discussed in lemma \ref{lem4n1}. Since the number of connections $\kappa_{j_mj_{m+1}}$
involved in the map $g$ is larger than the total number of connections of $Y$,
at least one connection occurs more than once.

Here we have two possibilities: either $g$ overlaps with a cyclic map
and we have $g_{s}=g_{d+s}$ for all $d< s< M$ or at least one of the equilibria
involved in $g$ has two incoming connections, say $\kappa_{j_mj_{m+1}}$ and
$\kappa_{j_sj_{s+1}}$. In the former case, the connection $\kappa_{ij}$
belongs to a minus-principal cycle, which contradicts the conditions in
the statement of this lemma. In the latter case due to lemmas \ref{lem4n1} and
\ref{lem4n2} the set $U_{Mk}$ is a thin cusp.
\qed

\begin{theorem}\label{th4n1}
Suppose that $Y$ is f.a.s. Generically, for sufficiently small $\delta>0$ almost all
$x\in\cB_{\delta}(Y)$ satisfy $\omega(x)=X$,
where $X\subset Y$ is a principal cycle.
\end{theorem}

\proof
To prove the theorem, it is sufficient to consider $x\in H_{ij}$.
If $\kappa_{ij}$ does not belong to a minus-principal cycle,
by lemma \ref{lem4n3} for almost all $x\in\cB_{\delta}(Y)\cap H_{ij}$
the trajectory $\Phi(x,t)$ for some $t>0$ (depending on $x$) belongs to one of
the sets $\{U_{M1},...,U_{MK}\}$ that are thin cusps. Hence, by corollary
\ref{principal} the trajectory
is attracted by a principal cycle as $t\to\infty$.

If $\kappa_{ij}$ belongs to a minus-principal cycle which is not
f.a.s., almost all trajectories $\Phi(x,t)$ with $x\in\cB_{\delta}(Y)\cap  H_{ij}$
for some $t>0$ escape from the $\delta$-neighbourhood of this cycle. They escape
along connections that do not belong to other minus-principal cycles
(two minus-principal cycles do not have common equilibria), hence the arguments
given above imply that they are attracted by principal cycles.

If a minus-principal cycle is f.a.s.\ (and thus principal), then the above arguments apply to those
trajectories that escape from it.
\qed

\begin{corollary}
If $Y$ is f.a.s., then its principal cycles that are e.a.s.\ are Milnor attractors.
\end{corollary}

\begin{corollary}\label{th41}
Consider a type A network $Y$ in $\R^4$. The network is f.a.s.\ if and only if
it contains at least one subcycle that is f.a.s.
\end{corollary}

Note that a corresponding result for networks of type Z in $\R^4$ follows from \cite[Theorem 3.4]{cl16}: two cycles in a network of type Z always have a connection in common, but it is not possible for a trajectory to switch from a neighbourhood of one cycle to the other and back. Therefore, whenever such a network is f.a.s., one of its subcycles must be f.a.s.

\section{Numerical examples: attracting heteroclinic cycles in
$(\D_2\rl\Z_4;\D_2\rl\Z_4)$-equivariant systems}\label{sec5}

In this section we construct two $(\D_2\rl\Z_4;\D_2\rl\Z_4)$-equivariant systems
with attracting heteroclinic cycles. According to theorem \ref{th32},
the graph associated with this group is of type IV. We aim at constructing
systems with type IVc networks, where in one system the four-equilibria subcycle
is stable, while in the other both two-equilibria subcycles are stable.
To achieve this we employ ideas from lemma 1 in \cite{clp17} that
we recall in subsection \ref{sec51}. The construction itself is discussed in
subsection \ref{sec52}.

\subsection{Construction of a $\Gamma$-equivariant system with a
desirable network}\label{sec51}

In \cite{clp17} we proved that a group $\Gamma$ with a set of invariant planes
satisfying certain conditions admits pseudo-simple heteroclinic cycles.
In the proof we explicitly built a $\Gamma$-equivariant
dynamical system $\dot {\bf x}={\bf f}({\bf x})$ possessing the desirable
cycle. As noted {\it ibid}, the proof can be generalised for the construction
of systems with other types of cycles and networks. Below we describe such
a construction for a given $\Gamma\subset \nSO(4)$ admitting a simple heteroclinic
network. Moreover, by varying constants employed in the construction,
we can change the eigenvalues of $\rd f(\xi_j)$, giving us control over the asymptotic stability of the cycles in the network.

We assume that the network under construction is given as a directed
graph, comprised of heteroclinic connections $\kappa_j:\xi_j\to\xi'_j$,
$1\le j\le m$. By $P_j=\Fix \Sigma_j$ we denote the isotropy plane
the connection belongs to, by $L_j$ and $L'_j$ the isotropy semiaxes in $P_j$
(since the network is simple, there are exactly two isotropy types of semiaxes
for each $P_j$), by $N_{\Gamma}(\Sigma_j)$,
$N_{\Gamma}(\Sigma_j)/\Sigma_j\cong\D_{K_j}$, the normalizer
of $\Sigma_j$ in $\Gamma$. If $P_i$ and $P_j$ intersect, then one of $L_i$ or
$L'_i$ coincides with one of $L_j$ or $L'_j$. The construction is done in three steps:

As a first step, for each $P_j$ we define a two-dimensional vector field
${\bf h}_j$, which in polar coordinates $(r,\theta)$ is:
\begin{equation}\label{fieh}
{\bf h}_j(r,\theta)=\left(r(1-r), \ \sin(K_j\theta)
(A_{j1}+A_{j2}\cos(K_j\theta))\right).
\end{equation}
We choose the angle of $L_j$ to be $\theta=0$, hence the angles of the nearest
$L'_j$ are $\pm\pi/K_j$. We assume that the sign of $A_{j1}+A_{j2}$ is positive and the sign
of $-A_{j1}+A_{j2}$ is negative. In \cite{clp17} we assumed that $A_{j2}=0$.
Here non-vanishing values of $A_{j2}$ are taken, because by varying this
parameter we change the eigenvalues of the linearisation $\rd f$ at $\xi_j$ and
$\xi'_j$, hence, change the stability properties of the cycles.
For the flow of $(\dot r,\dot\theta)={\bf h}_j(r,\theta)$ each of the axes
$\theta=2k\pi/K_j$ ($L_j$) and $\theta=(2k+1)\pi/K_j$ ($L'_j$) is invariant and
has an equilibrium at $r=1$ which is attracting along the direction of $r$.
Moreover, there are heteroclinic connections between equilibria on
neighbouring axes, going from $\xi_j\in L_j$ to $\xi_j'\in L_j'$.

As a second step we extend the vector fields ${\bf h}_j$ to
${\bf g}_j: \R^4 \to \R^4$ as follows: Denote by $\pi_j$ and $\pi^{\perp}_j$
the projections onto the plane $P_j$ and its orthogonal complement in $\R^4$,
respectively. We set
\begin{equation}\label{gj}
\pi_j{\bf g}_j({\bf x})=
{{\bf h}_j(\pi_j{\bf x})\over 1+B|\pi^{\perp}_j{\bf x}|^2},\quad
\pi^{\perp}_j{\bf g}_j({\bf x})=0,
\end{equation}
with a positive constant $B$, which should be taken sufficiently large.

Finally, we define the vector field ${\bf f}:\R^4\to\R^4$ as
\begin{equation}\label{f-equation}
{\bf f}({\bf x})=\sum\limits_{j=1}^{m}
\sum\limits_{\gamma_{ij}\in \mathcal{G}_j} \gamma_{ij}{\bf g}_j(\gamma_{ij}^{-1}{\bf x}),
\end{equation}
where $\mathcal{G}_j=\Gamma/N_{\Gamma}(\Sigma_j)$.

One can prove then that the system
\begin{equation}\label{f-system}
\dot{\bf x}={\bf f}({\bf x})
\end{equation}
possesses steady states $\xi_j\in L_j$, for all $L_j$ involved in the network,
and heteroclinic connections $\xi_j\to\xi_j'\subset P_j$.
By construction the system (\ref{f-system}) is $\Gamma$-equivariant,
which implies the invariance of all axes $L_j$ and planes $P_j$. The eigenvalues
of $\rd f(\xi_j)$ and $\rd f(\xi_j')$ associated with eigenvectors in $P_j$
for large values of $B$ are approximated by
$A_{j1}+A_{j2}>0$ and $-A_{j1}+A_{j2}<0$, respectively.

\subsection{$(\D_2\rl\Z_4;\D_2\rl\Z_4)$-equivariant systems with
attracting heteroclinic cycles}\label{sec52}

Using the results of the previous subsection, we construct two
$(\D_2\rl\Z_4;\D_2\rl\Z_4)$-equivariant dynamical systems possessing a type IVc network,
such that in the first system the four-equilibria cycle is attracting and in the
second both two-equilibria cycles are attracting.

The group $\Gamma=(\D_2\rl\Z_4;\D_2\rl\Z_4)$ is comprised of the following
elements:
\begin{eqnarray*}\label{listex1}
\kappa_1(\pm)=\pm((1,0,0,0);(1,0,0,0))&,&\kappa_2(\pm)=\pm((1,0,0,0);(0,0,0,1)),\\
\kappa_3(\pm)=\pm((0,0,0,1);(1,0,0,0))&,&\kappa_4(\pm)=\pm((0,0,0,1);(0,0,0,1)),\\
\kappa_5(\pm)=\pm((0,1,0,0);(0,1,0,0))&,&\kappa_6(\pm)=\pm((0,0,1,0);(0,0,0,1)),\\
\kappa_7(\pm)=\pm((0,1,0,0);(0,0,1,0))&,&\kappa_8(\pm)=\pm((0,0,1,0);(0,1,0,0)).
\end{eqnarray*}
The group has six isotropy types of subgroups $\Sigma$ satisfying
$\dim\Fix\Sigma=2$, the respective fixed-point subspaces are:
\begin{equation}\label{lissub}
\renewcommand{\arraystretch}{1.5}
\begin{array}{l|c|l}
\hbox{Plane}&\hbox{Fixed by}&\hbox{Coordinates}\\
\hline
P_1&\kappa_4(+)&(x_1,0,0,x_4)\\
P_2&\kappa_4(-)&(0,x_2,x_3,0)\\
P_3&\kappa_5(\pm)&(x_1,x_2,0,0),\ (0,0,x_3,x_4)\\
P_4&\kappa_6(\pm)&(x_1,0,x_3,0),\ (0,x_2,0,x_4)\\
P_5&\kappa_7(\pm)&(x_1,x_2,x_2,-x_1),\ (x_1,x_2,-x_2,x_1)\\
P_6&\kappa_8(\pm)&(x_1,x_2,x_2,x_1),\ (x_1,x_2,-x_2,-x_1).
\end{array}
\end{equation}
The invariant axes are:
\begin{equation}\label{axsub}
\renewcommand{\arraystretch}{1.5}
\begin{array}{l|c|l}
\hbox{Axis}&\hbox{Intersection of}&\hbox{Coordinates}\\
\hline
L_1&P_1,\ P_3,\ P_4&(x_1,0,0,0),\ (0,0,0,x_4)\\
L_2&P_2,\ P_3,\ P_4&(0,x_2,0,0),\ (0,0,x_3,0)\\
L_3&P_1,\ P_5,\ P_6&(x_1,0,0,x_1),\ (x_1,0,0,-x_1)\\
L_4&P_2,\ P_5,\ P_6&(0,x_2,x_2,0),\ (0,x_2,-x_2,0)
\end{array}
\end{equation}
Since $-I\in\Gamma$, two semi-axes of any axis are of the same isotropy type.

\begin{figure}

\vspace*{-12mm}
\hspace*{-20mm}\includegraphics[width=12cm]{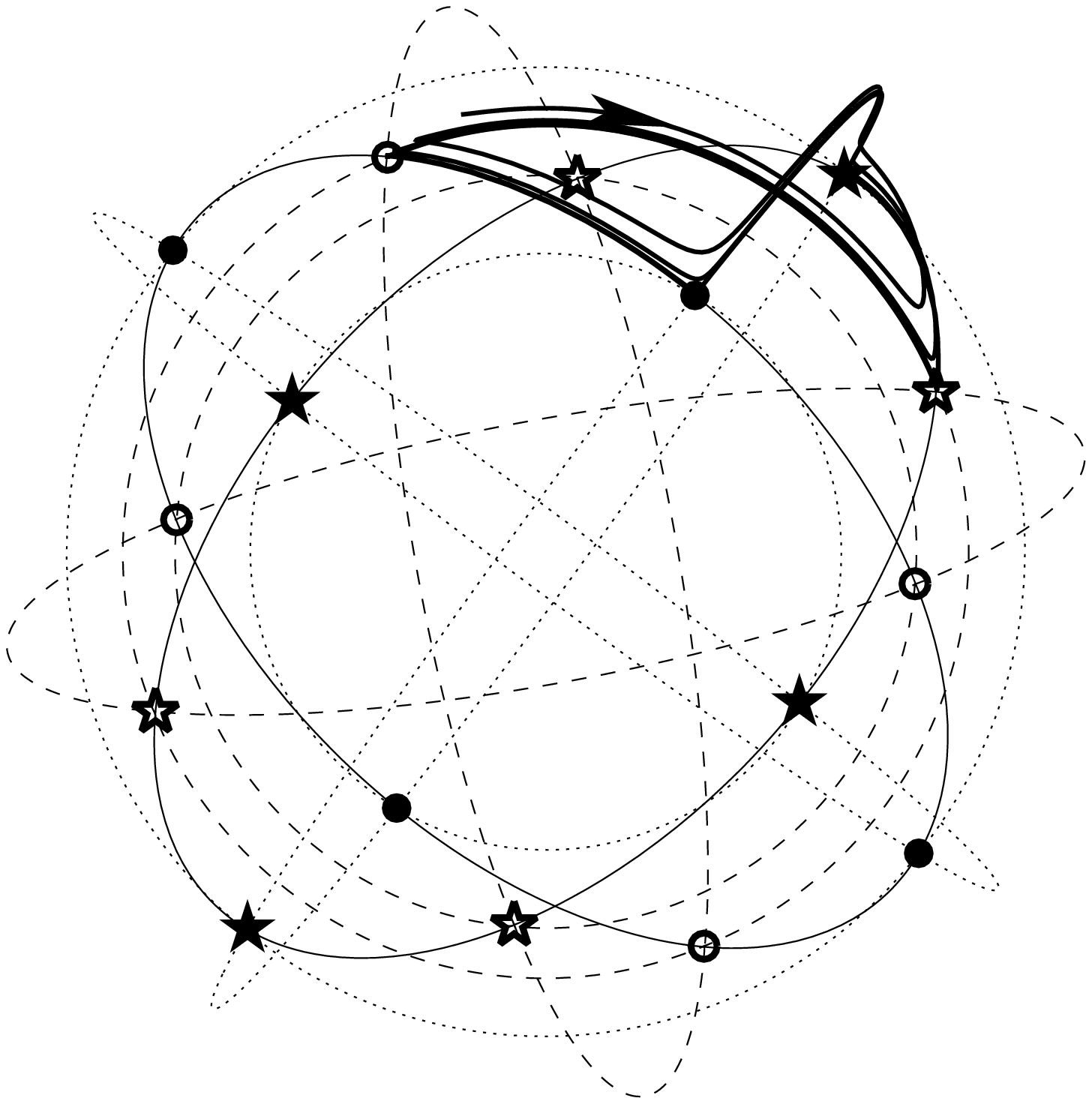}

\vspace*{-40mm}
\hspace*{100mm}{\Large (a)}

\vspace*{12mm}
\hspace*{50mm}\includegraphics[width=12cm]{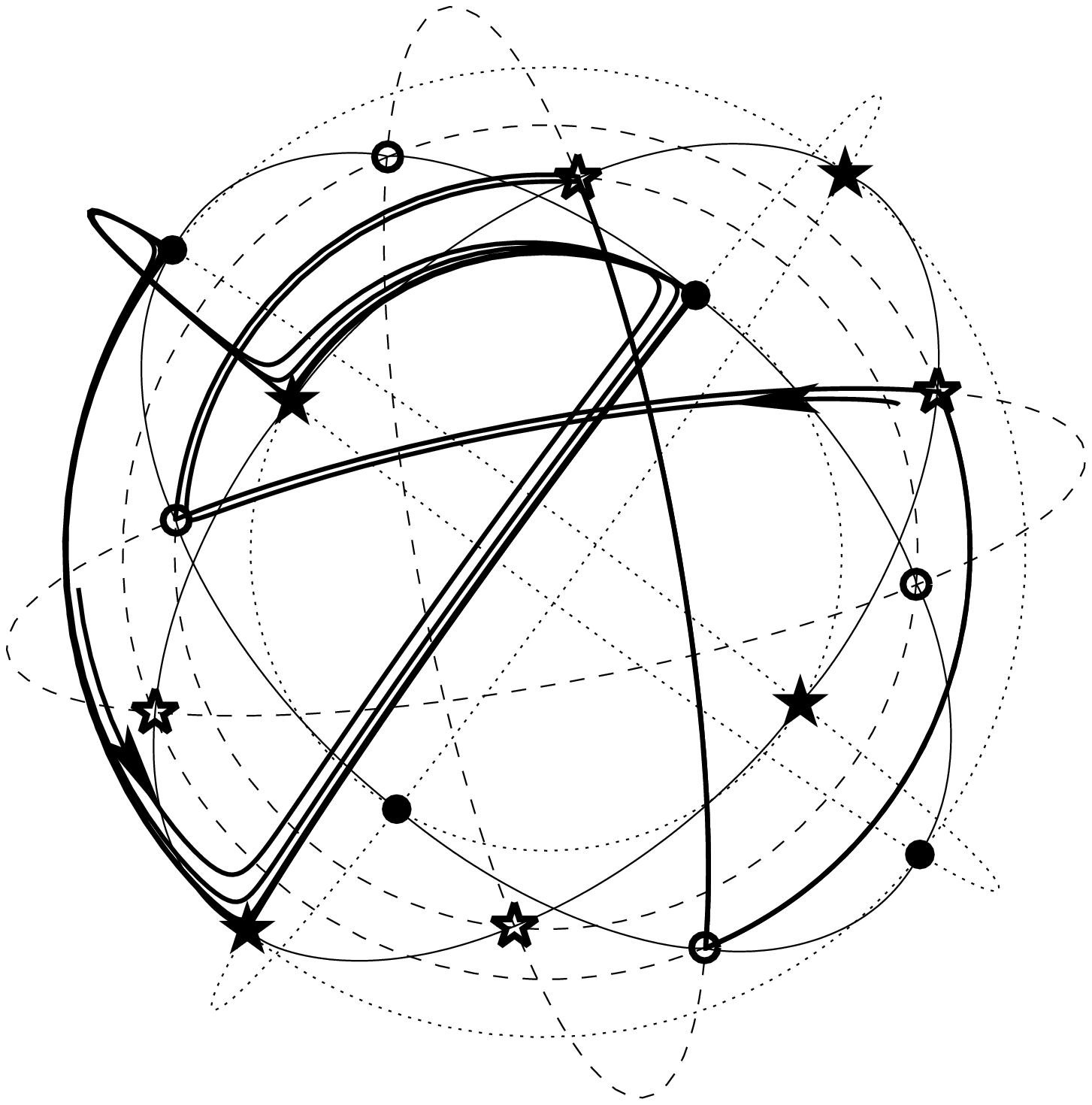}

\vspace*{-30mm}
\hspace*{60mm}{\Large (b)}

\vspace*{15mm}
\noindent
\caption{Projection of the heteroclinic connections in the planes $P_1$ and
$P_2$ (solid lines), $P_3$ and $P_4$ (dashed lines), $P_5$ and $P_6$
(dotted lines) and trajectories (bold lines) approaching the four-equilibria
cycle (a) and the two-equilibria cycles (b) into the plane
$<{\bf v}_1,{\bf v}_2>$, where ${\bf v}_1=(4,2,4,1.5)$ and
${\bf v}_2=(2,4,-1.5,4)$. The steady states $\xi_1$ and $\xi_2$ are denoted by
filled circles and stars, the steady states $\xi_3$ and $\xi_4$ by hollow
circles and stars, respectively.}
\label{fig51}\end{figure}

We have $K_1=K_2=4$ and $K_j=2$ for $3\le j\le 6$, therefore in agreement with
(\ref{fieh}) set
\begin{eqnarray*}
{\bf h}_j(r,\theta)&=&\left(r(1-r), \ \sin4\theta
(A_{j1}+A_{j2}\cos4\theta)\right)\hbox{ for }j=1,2,\\
{\bf h}_j(r,\theta)&=&\left(r(1-r), \ \sin2\theta
(A_{j1}+A_{j2}\cos2\theta)\right)\hbox{ for }3\le j\le 6.
\end{eqnarray*}
For the plots shown in figure \ref{fig51} the coefficients $A_{j1}$ and $A_{j2}$ are:
\begin{align*}
(a)&:~~A_{11}=A_{21}=25,\ A_{12}=A_{22}=-5,\
A_{j1}=15,\ A_{j2}=-5,\hbox{ for }3\le j\le 6,\\
(b)&:~~A_{11}=A_{21}=2,\ A_{12}=A_{22}=-1,\
A_{j1}=25,\ A_{j2}=-5,\hbox{ for }3\le j\le 6,
\end{align*}
and in both cases $B=100$. Recall that a type A heteroclinic cycle is e.a.s.\ if
\begin{equation}\label{condst}
\prod_{1\le j\le J}{|c_j|\over e_j}>1\hbox{ and }e_j>t_j\hbox{ for all }1\le j\le J,
\end{equation}
where we use the usual convention of writing $-c_j, e_j$ and $t_j$ for the contracting, expanding and transverse eigenvalues, respectively. In case (a) the eigenvalues of $\rd f(\xi_j)$ for the four-equilibria cycle are:
\begin{eqnarray*}
e_1=e_3\approx10,\ c_1=c_3\approx-30,\ t_1=t_3\approx-20,\
e_2=e_4\approx20,\ c_2=c_4\approx-20,\ t_2=t_4\approx10,\
\end{eqnarray*}
therefore (\ref{condst}) implies that the cycle is e.a.s.
In case (b) for both two-equilibria cycles the eigenvalues are
$$e_1=e_2\approx20,\ c_1=c_2\approx-30,\ t_1=-3,\ t_2=1,$$
hence (\ref{condst}) implies that these cycles are e.a.s.
For both cases trajectories approaching attracting the heteroclinic cycles are shown
on figure \ref{fig51}.

\section{Conclusion}\label{sec6}
We have contributed to the systematic study of heteroclinic dynamics in low dimensions in two ways: (i) we derived a complete list of maximal simple heteroclinic networks in $\R^4$; and (ii) we proved conditions for fragmentary and essential asymptotic stability of type A heteroclinic networks in $\R^4$, i.e.\ admitted by a group $\Gamma\subset \nSO(4)$. Along the way we introduced the concept of a graph associated with a group, presented a list of simple graphs for subgroups of $\nO(4)$ and for each graph the subgroups that the graph is associated with.

Possible continuations of this work include a similar study for pseudo-simple
heteroclinic networks in $\R^4$, extending the work on pseudo-simple cycles in \cite{clp17}.
%In principle, the same classification of graphs and eventually heteroclinic cycles/networks associated to groups $\Gamma \subset \nO(n)$ can be attempted for higher dimensions.
The employed approach is likely to be applicable for the classification of graphs and
eventually heteroclinic  cycles/networks associated with groups
$\Gamma \subset \nO(n)$ for $n>4$. By contrast, the proofs of the stability results
rely on properties of thin and thick cusps, for which the generalisation
to higher dimensions is not obvious.

\medbreak
{\bf Acknowledgements}: The second author gratefully acknowledges support through the project 57338573 PPP Portugal 2017 by the German Academic Exchange Service (DAAD), sponsored by the Federal Ministry of Education and Research (BMBF).

\end{document}